\RequirePackage[2020-02-02]{latexrelease}
\documentclass[11pt]{amsart}
\usepackage{amsfonts}
\usepackage{amsmath}
\usepackage{amsthm}
\usepackage{amssymb}
\usepackage{wasysym}
\usepackage{mathrsfs}
\usepackage[numbers]{natbib}
\usepackage[fit]{truncate}
\usepackage{hyperref}
\usepackage{theoremref}

\usepackage{xcolor}

\usepackage{tikz-cd}

\usepackage{thmtools}
\usepackage{thm-restate}

\usepackage{fmtcount}

\usepackage{enumitem}

\usepackage[UKenglish]{babel}
\usepackage[UKenglish]{isodate}

\usepackage{verbatim}
\usepackage[utf8]{inputenc}

\usepackage{dsfont}

\usepackage{etoolbox}
\makeatletter
\patchcmd{\@maketitle}
{\ifx\@empty\@dedicatory}
{\ifx\@empty\@date \else {\vskip3ex \centering\footnotesize\@date\par\vskip1ex}\fi
	\ifx\@empty\@dedicatory}
{}{}
\patchcmd{\@adminfootnotes}
{\ifx\@empty\@date\else \@footnotetext{\@setdate}\fi}
{}{}{}
\makeatother

\usepackage{catoptions}
\makeatletter

\def\Autoref#1{%
	\begingroup
	\edef\reserved@a{\cpttrimspaces{#1}}%
	\ifcsndefTF{r@#1}{%
		\xaftercsname{\expandafter\testreftype\@fourthoffive}
		{r@\reserved@a}.\\{#1}%
	}{%
		\ref{#1}%
	}%
	\endgroup
}
\def\testreftype#1.#2\\#3{%
	\ifcsndefTF{#1autorefname}{%
		\def\reserved@a##1##2\@nil{%
			\uppercase{\def\ref@name{##1}}%
			\csn@edef{#1autorefname}{\ref@name##2}%
			\autoref{#3}%
		}%
		\reserved@a#1\@nil
	}{%
		\autoref{#3}%
	}%
}
\makeatother

\numberwithin{equation}{section}

 \newtheoremstyle{dotless}{}{}{\itshape}{}{\bfseries}{}{ }{}

\theoremstyle{dotless}

\MakeRobust{\ref}% avoid expanding it when in a textual label

\makeatletter
\newcommand{\labeltext}[2]{%
	\@bsphack
	\csname phantomsection\endcsname % in case hyperref is used
	\def\@currentlabel{#1}{\label{#2}}%
	\@esphack
}
\makeatother

\makeatletter
\def\blfootnote{\gdef\@thefnmark{}\@footnotetext}
\makeatother

\makeatletter
\def\moverlay{\mathpalette\mov@rlay}
\def\mov@rlay#1#2{\leavevmode\vtop{%
		\baselineskip\z@skip \lineskiplimit-\maxdimen
		\ialign{\hfil$\m@th#1##$\hfil\cr#2\crcr}}}
\newcommand{\charfusion}[3][\mathord]{
	#1{\ifx#1\mathop\vphantom{#2}\fi
		\mathpalette\mov@rlay{#2\cr#3}
	}
	\ifx#1\mathop\expandafter\displaylimits\fi}
\makeatother

\newcommand{\bigcupdot}{\charfusion[\mathop]{\bigcup}{\cdot}}

\newcommand{\N}[0]{\mathbb{N}}
\newcommand{\Z}[0]{\mathbb{Z}}
\newcommand{\Q}[0]{\mathbb{Q}}
\newcommand{\R}[0]{\mathbb{R}}

\newcommand{\OO}[0]{R}
\newcommand{\II}[0]{I}

\newcommand{\bA}[0]{\mathbf{A}}

\newcommand{\TL}[0]{T_{\mathrm{L}}}
\newcommand{\supp}[0]{\mathrm{supp}}

\newcommand{\U}{\mathcal{U}}
\newcommand{\brackets}[1]{\left( #1 \right)}

\newcommand{\ip}[1]{\left\lfloor #1\right\rfloor}

\newcommand{\Lor}[0]{\mathcal{L}_{\mathrm{or}}}

\newcommand{\Lotf}[0]{\mathcal{L}_{\mathrm{otf}}}

\newcommand{\Totf}[0]{T_{\mathrm{otf}}}

\newcommand{\ol}[1]{\overline{#1}}

\newcommand\restr[2]{{
		\left.\kern-\nulldelimiterspace 
		#1
		\vphantom{\big|} 
		\right|_{#2}
}}

\newcommand{\one}[0]{\mathds{1}}
\newcommand{\cal}[1]{\mathcal{#1}}
\newcommand{\PP}[0]{\mathbf{P}}

\DeclareMathOperator{\Th}{Th}

\usepackage{theoremref}

\theoremstyle{definition}
\newtheorem{definition}{Def\-i\-ni\-tion}[section]

\newtheorem{example}[definition]{Example}
\newtheorem{remark}[definition]{Remark}
\newtheorem{construction}[definition]{Construction}
\theoremstyle{plain}
\newtheorem{theorem}[definition]{Theorem}

\newtheorem{proposition}[definition]{Prop\-o\-si\-tion}
\newtheorem{lemma}[definition]{Lem\-ma}
\newtheorem{corollary}[definition]{Corollary}
\newtheorem{question}[definition]{Question}

\begin{document}

	\begin{abstract}
		We develop a first-order theory of ordered transexponential fields in the language $\{+,\cdot,0,1,<,e,T\}$, where $e$ and $T$ stand for unary function symbols. While the archimedean models of this theory are readily described, the study of the non-archimedean models leads to a systematic examination of the induced structure on the residue field and the value group under the natural valuation. We establish necessary and sufficient conditions on the value group of an ordered exponential field $(K,e)$ to admit a transexponential function $T$ compatible with $e$. Moreover, we give a full characterisation of all \emph{countable} ordered transexponential fields in terms of their valuation theoretic invariants.
	\end{abstract}

	\title{Ordered transexponential fields}
	
	\author[L.~S.~Krapp]{Lothar Sebastian Krapp}%\thanks{Supported by...}
	\author[S.~Kuhlmann]{Salma Kuhlmann}
	
	\address{Fachbereich Mathematik und Statistik, Universität Konstanz, 78457 Konstanz, Germany}
	\email{sebastian.krapp@uni-konstanz.de}
	
	\address{Fachbereich Mathematik und Statistik, Universität Konstanz, 78457 Konstanz, Germany}
	\email{salma.kuhlmann@uni-konstanz.de}

	\date{\today}
	
	\cleanlookdateon
	
	\maketitle
	
	%\tableofcontents
	
	\blfootnote{\textup{2020} \textit{Mathematics Subject Classification}: 
		12J25, 03C64 (primary); 12J15, 26A12, 06F20, 12L12 (secondary)
		}
	
\section{Introduction}

	During the 1980s and 1990s, the algebraic and model theoretic properties of the real exponential field $\R_{\exp}=(\R,+,\cdot,0,1,<,\exp)$, where $\exp$ denotes the standard exponential $x\mapsto \mathrm{e}^x$, were studied extensively. This area of research, going back to Tarski's highly influential work \cite{tarski}, was mainly motivated by establishing analogies between the complete first-order theory of $\R_{\exp}$ and the complete first-order theory of the ordered field of real numbers $(\R,+,\cdot,0,1,<)$, now known as the first-order theory of real closed fields. A major breakthrough was made by Wilkie in \cite{wilkie}, who proved that $\R_{\exp}$ is model complete and thus o-minimal (which was also established by van den Dries, Macintyre and Marker in  \cite{dries2}). 
	Studying the growth properties of definable functions in o-minimal expansions of the ordered field of real numbers \cite{miller2,miller3} and in o-minimal expansions of any ordered field \cite{miller}, Miller established the following remarkable growth dichotomy: an o-minimal expansion of an ordered field is either power bounded or admits a definable exponential function. 
	Going one step further in the hierarchy of growth, Miller's dichotomy result naturally led to the question whether there exist o-minimal expansions of ordered fields that are not exponentially bounded.
	Recent research activity in this area is therefore motivated by the search for \emph{either} an o-minimal expansion of an ordered exponential field by a transexponential function that eventually exceeds any iterate of the exponential \emph{or}, contrarily, for a proof that any o-minimal expansion of an ordered field is already exponentially bounded. (See also \cite[Section~1.1]{padgett} for further background.)
	
	In this paper, we initiate a valuation theoretic approach towards the study of ordered exponential fields endowed with transexponential functions in order to gain a better understanding of transexponential growth. Our approach is somewhat similar to the one chosen by the second author in \cite{kuhlmann} towards establishing a theory of ordered exponential fields. Namely, we fix a simple axiomatisation of the theory of ordered transexponential fields that we consider and then study the valuation theoretic invariants of the non-archimedean models of this theory. In particular, we examine the structure that a transexponential induces on the residue field and the value group under the natural valuation. Our notion of transexponentials is inspired by Boshernitzan~\cite{boshernitzan}, that is, given an exponential $e$, a transexponential $T$ has the property $T(x+1)=e(T(x))$ for any non-negative $x$.
	
	The structure and the main achievements of this paper are as follows. After establishing some general preliminaries on ordered structures, valuations and exponentials in \Autoref{sec:prelim}, we introduce the theory $\Totf$ of ordered transexponential fields in \Autoref{sec:theory}. We provide an axiomatisation (\Autoref{def:totf}) and introduce the growth property (\Autoref{def:grwoth}) for transexponentials. Moreover, we prove some basic analytic properties of transexponentials used in \Autoref{thm:notomin} to show that $\Totf$ does not admit any o-minimal models.
	In \Autoref{sec:res}, we consider the restriction of transexponentials to the valuation ring and obtain in \Autoref{thm:restrans} that its induced map on the (archimedean) residue field is also a transexponential. Moreover, \Autoref{thm:arch} shows that any archimedean ordered transexponential field is already a substructure of the real transexponential field (see \Autoref{ex:realtrans}). In \Autoref{sec:pos}, we consider the restriction of transexponentials to the ordered set of positive infinite elements of the field. 
	Under sufficient conditions on the principal exponential rank of an ordered exponential field $(K,e)$, we
	 present in \Autoref{cons:trans} a general procedure to obtain \emph{all} transexponentials on $K$ compatible with $e$ (see \Autoref{thm:alltransexp}). \Autoref{sec:val} is ded\-i\-cat\-ed to the study of the structure transexponentials induce on the value group. In \Autoref{prop:condextransexp2}, we establish necessary and sufficient conditions on the value group of an ordered exponential field (respectively its rank) in order that it admits a transexponential. Finally, \Autoref{sec:countable} studies countable non-archimedean ordered transexponential fields. We construct countable ordered exponential fields admitting no compatible transexponentials (see \Autoref{rmk:constr}) as well as transexponentials without the growth property (see \Autoref{ex:nogrowth}), a key step for which is the construction of countable ordered exponential fields with prescribed principal exponential rank in \Autoref{con:exprank}.
	While \Autoref{sec:res} to \Autoref{sec:val} present several `going down' approaches -- that is, general (non-archimedean) models of $\Totf$ are considered and then the structure induced by the transexponential on the residue field and the value group are studied --, \Autoref{sec:countable}  initiates a `going up' approach by starting with a given archimedean field and ordered abelian group to construct \emph{countable} models of $\Totf$ (\Autoref{thm:nonarchexpand}). Other constructions in the \emph{uncountable} case via power series will be the subject of a future publication \cite{krappkuhlmannwip}.

\section{Preliminaries}\label{sec:prelim}

	We follow closely the notation and terminology of \cite{kuhlmann}. For the convenience of the reader, we introduce the relevant background on ordered structures, valuation theory and exponentials below.
	
	\subsection{Ordered structures}
	We denote by $\N$ the set of natural numbers without $0$.  Given an (additive) ordered abelian group $H$, for any $h\in H$ we denote by $|h|$ the absolute value of $h$, i.e.\ $|h|=\max\{h,-h\}$.
	
	The language of ordered rings $\{+,\cdot,0,1,<\}$ is denoted by $\Lor$.
	If no confusion is likely to arise, the $\Lor$-structure of an ordered field $(K,+,\cdot,0,1,<)$ is simply abbreviated by $K$. Likewise, for any expansion of $(K,+,\cdot,0,1,<)$ we often omit the symbols $\{+,\cdot,0,1,<\}$. We use similar conventions for ordered abelian groups, i.e.\ we write $H$ instead of $(H,+,0,<)$. 
	When we say `definable' in the model theoretic sense, we always mean `definable with parameters'. 
	\textbf{Throughout this work, $K$ denotes an ordered field.}
	
	Let $S$ denote a set totally ordered by $<$. We use standard interval notation in $S$, e.g.\ $[a,b]_S$ for the interval $\{x\in S\mid a\leq x\leq b\}$ where $a,b\in S$. Moreover, we use superscript notation like $S^{>a}$ to denote infinite intervals, e.g.\  $S^{>a}=(a,\infty)_S=\{x\in S\mid a<x\}$. 
	For any monotonically increasing map $\varphi\colon S\to S$ (i.e.\ $\varphi(a)\leq \varphi(b)$ if $a<b$) and any $n\in \N\cup\{0\}$, we denote by $\varphi^n$ the \textbf{$n$-th iterate} of $\varphi$, i.e.\ $\varphi^0(a)=a$ and $\varphi^{n+1}(a)=\varphi(\varphi^n(a))$ for any $a\in S$. 
	If, moreover, $\varphi$ is bijective, then $\varphi^{-n}$ denotes the $n$-th iterate of the inverse $\varphi^{-1}$ of $\varphi$.
	Two elements $a,b\in S$ are called \textbf{$\varphi$-equivalent} if there exists some $n\in \N$ such that $a\leq \varphi^n(b)$ and $b\leq \varphi^n(a)$. We then write $a\sim_\varphi b$ and denote the equivalence class of $a$ by $[a]_{\varphi}$. Since $[a]_{\varphi}$ is a convex subset of $S$, we obtain an ordering on $S/\!\sim_{\varphi}$ given by $[a]_{\varphi}<[c]_{\varphi}$ if and only if $a<c$ and $a\not\sim_{\varphi}c$, for any $a,c\in S$. See \cite[Remark~3.20]{kuhlmann} for further details.
		
	Recall that due to Hölder's Theorem (cf.\ e.g.\ \cite[Theorem~1.1]{kuhlmann}), any archimedean ordered field $K$ embeds into $\R$ via a unique $\Lor$-embedding. We often identify $K$ with the corresponding subfield of $\R$.
	
	Given a non-empty ordered set $\Delta$, we denote by $$\coprod_\Delta K$$
	the \textbf{Hahn sum} over $\Delta$ and $K$. It consists of all maps $s\colon \Delta \to K$ whose support $\supp(s)=\{a\in \Delta \mid s(a)\neq 0\}$ is finite. We endow $\coprod_\Delta K$ with pointwise addition and the lexicographic order $s>0 \Leftrightarrow s(\min\supp(s))>0$ for any $s\neq 0$, making it an ordered abelian group. For a given $\delta\in \Delta$ we denote by $\one_\delta\colon \Delta \to K$ the indicator function mapping $\delta$ to $1$ and everything else to $0$. We may then express an element $s\in \coprod_\Delta K$ by $s=\sum_{\delta\in \Delta}s(\delta)\one_\delta$. 
	
	\subsection{Valuations}
	
	For a valuation $w$ on an ordered field $K$, we denote its value group by $w(K)$ (written additively), its valuation ring $\{x\in K\mid w(x)\geq 0\}$ by $R_w$, its valuation ideal $\{x\in K\mid w(x)>0\}$ by $I_w$, its (multiplicative) group of units $\{x\in K\mid w(x)=0\}$ by $\mathcal{U}_w$ and its (multiplicative) group of $1$-units $\{1+x\mid x\in I_w\}$ by $1+I_w$.
	
	The \textbf{natural valuation} $v$ on $K$ is the finest non-trivial convex valuation on $K$, i.e.\ the valuation on $K$ whose valuation ring $R_v$ is the convex hull of $\Z$ in $K$. Its residue field $R_v/I_v$ is simply denoted by $\ol{K}$, and for any $a\in R_v$, we also denote $a+I_v\in \ol{K}$ by $\ol{a}$. Note that $\ol{K}$ is an ordered field whose ordering is given by $\ol{a}<\ol{b} \Leftrightarrow (a<b\wedge \ol{a}\neq \ol{b})$ for any $a,b\in R_v$. The value group $v(K)$ is usually denoted by $G$.
	
	Elements of $I_v$ are called \textbf{infinitesimals}, of $R_v$ \textbf{finite elements} and of $\mathbf{P}_K=\{x\in K^{>0}\mid v(x)<0\}$ \textbf{positive infinite elements}.
	We obtain the following picture:

	\begin{center}
		\setlength{\unitlength}{0.002\textwidth}
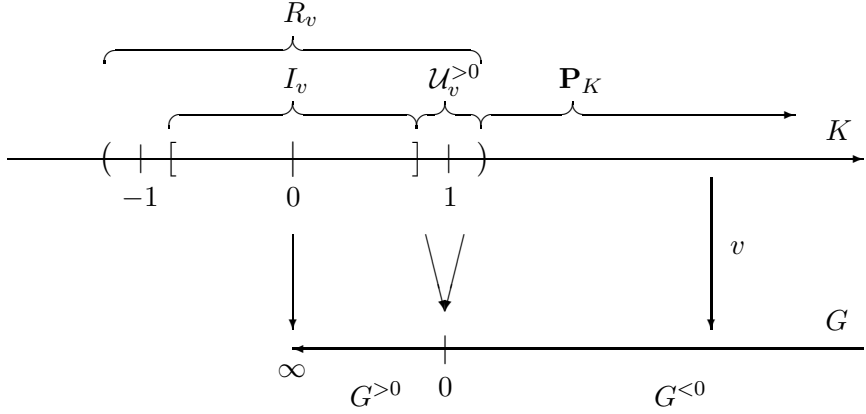
		\begin{picture}(460,250)(-150,-50)
			\put(-102,96){{$\left(\hspace{20pt}\left[\hspace{42pt}\left|
					\rule{0pt}{9pt} \right.\hspace{42pt}\right]\hspace{20pt}\right)$}}
			\put(-150,100){\vector(1,0){450}}
			\put(300,0){\vector(-1,0){300}}
			\put(220,90){\vector(0,-1){80}}
			\put(230,50){{$v$}}
			\put(0,60){\vector(0,-1){50}}
			\put(70,60){\vector(1,-4){10}}
			\put(90,60){\vector(-1,-4){10}}
			\put(-3.5,75){{$0$}}			
			\put(280,110){{$K$}}
			\put(280,10){{$G$}}
			\put(-8,-15){{$\infty$}}
			\put(30,-30){{$G^{>0}$}}
			\put(190,-30){{$G^{<0}$}}
			\put(80,96){{$|$}}
			\put(79,75){{$1$}}
			\put(-82,96){{$|$}}
			\put(-90,75){{$-1$}}
			\put(78,-4){{$|$}}
			\put(76.5,-25){{$0$}}
			\put(-92,157){\line(1,0){85}}
			\put(7,157){\line(1,0){85}}
			\put(-92,150){\oval(14,14)[tl]}
			\put(92,150){\oval(14,14)[tr]}
			\put(-7,164){\oval(14,14)[br]}
			\put(7,164){\oval(14,14)[bl]}
			\put(-5,172){{$R_v$}}
			\put(-58,123){\line(1,0){51}}
			\put(7,123){\line(1,0){51}}
			\put(-58,116){\oval(14,14)[tl]}
			\put(58,116){\oval(14,14)[tr]}
			\put(-7,130){\oval(14,14)[br]}
			\put(7,130){\oval(14,14)[bl]}
			\put(-5,137){{$I_v$}}
			\put(72,123){\line(1,0){1}}
			\put(87,123){\line(1,0){5}}
			\put(72,116){\oval(14,14)[tl]}
			\put(92,116){\oval(14,14)[tr]}
			\put(73,130){\oval(14,14)[br]}
			\put(87,130){\oval(14,14)[bl]}
			\put(72,137){{${\cal U}_v^{>0}$}}
			\put(106,123){\line(1,0){34}}
			\put(154,123){\vector(1,0){110}}
			\put(106,116){\oval(14,14)[tl]}
			\put(140,130){\oval(14,14)[br]}
			\put(154,130){\oval(14,14)[bl]}
			\put(140,137){{{\bf P}$_K$}}
		\end{picture}
	\end{center}

	Given two ordered abelian groups $H_1$ and $H_2$, we denote by $H_1\amalg H_2$ their sum $H_1\oplus H_2$ equipped with the lexicographic order. Using this notation, we fix for $(K,+,0,<)$ an additive decomposition $$K = \bA \amalg R_v,$$
	where $\bA$ is a complement to the valuation ring $R_v$ (cf.\ \cite[Theorem~1.4]{kuhlmann}). Note that $\bA$ is unique up to order-preserving isomorphism and $v(\bA)=G^{<0}$. 
	
	We denote the \textbf{natural valuation} on $G$ by $v_G$. It maps each element $g\in G$ to its archimedean equivalence class $[g]_+=\{h\in G\mid |g|\leq n|h|\leq n^2|g|\text{ for some } n\in \N\}$. %Note that the archimedean equivalence class of some $g\in G$ is given by $[g]_\varphi\cup [-g]_\varphi$ for the map $\varphi\colon G\to G, g\mapsto 2g$. 
	The \textbf{rank} of $G$, i.e.\ the ordered value set $v_G(G)=\{[g]_+\mid g\in G\setminus\{0\}\}$, is denoted by $\Gamma$. Its ordering is given by $[g]_+<[h]_+\Leftrightarrow (|g|>|h|\wedge [g]_+\neq[h]_+)$ for any $g,h\in G\setminus \{0\}$. 
		
	\subsection{Exponentials}
	
	A unary function $e\colon K \to K^{>0}$ is called an \textbf{exponential} on $K$ if it is an order-preserving isomorphism from the ordered additive group $(K,+,0,<)$ of $K$ to the ordered multiplicative group $(K^{>0},\cdot,1,\allowbreak<)$ of positive elements of $K$. We then call $(K,e)$ an \textbf{ordered exponential field}. 
	The inverse of $e$ is called a \textbf{logarithm} and denoted by $\ell$. 
	A subfield $F$ of $K$ is called \textbf{exponentially closed} in $(K,e)$ if the restriction $e|_F$ is an exponential on $F$. 
	Throughout, we focus on exponentials satisfying the growth axiom scheme (GA) and the Taylor axiom (T$_1$):
	\begin{itemize}
		\item[(GA)] for any $n\in \N$ and for any $a\in K$, if $a\geq n^2$, then $e(a)>a^n$;
		
		\item[(T$_1$)] for any $a\in K$ with $0<|a|\leq 1$, we have $|e(a)-(1+a)|<a^2$.
	\end{itemize}
	Such an exponential is called a \textbf{(GAT$_{\mathbf{1}}$)-exponential}.
	Likewise, an \textbf{ordered (GAT$_{\mathbf{1}}$)-exponential field} is an ordered exponential field $(K,e)$ whose exponential satisfies (GA) and (T$_1$).
	
	An exponential $e$ is called \textbf{$v$-compatible} if its \textbf{residue exponential} $\ol{e}\colon \ol{K}\to\ol{K}^{>0}$ given by $\ol{e}(\ol{a})=\ol{e(a)}$ for any $a\in R_v$ is a well-defined exponential on $\ol{K}$ (see \cite[Lemma~1.17]{kuhlmann} for equivalent definitions). When we consider differentiability of $e$ (or other unary functions on $K$), we mean it in a first-order sense, that is, for given $a,b\in K$, the equation $e'(a)=b$ expresses that for any $\varepsilon \in K^{>0}$, there exists $\delta\in K^{>0}$ such that for any $h\in K$ with $0<|h|<\delta$,
	$$\left|\frac{e(a+h)-e(a)}{h}-b\right|<\varepsilon.$$
	Due to \cite[Theorem~3.40]{krappthesis}, we have the following.
	
	\begin{lemma}\thlabel{lem:t1exp}
		Any exponential $e$ on $K$ for which (T$_1$) holds is $v$-compatible and satisfies the differential equation $e'=e$.
	\end{lemma}

	In \cite[Theorem~3.3]{krappthesis}, a characterisation of all archimedean ordered exponential fields is given. Combining this with \Autoref{lem:t1exp}, we obtain the following.

	\begin{lemma}\thlabel{lem:gatgat}
		Let $(K,e)$ be an ordered (GAT$_{{1}}$)-exponential field. Then $(\ol{K},\ol{e})\subseteq (\R,\exp)$. In particular,
		also $\ol{e}$ is a (GAT$_{{1}}$)-exponential on $\ol{K}$.
	\end{lemma}

	\begin{proof}
		First note that $\ol{e}$ satisfies $|\ol{e}(x)-(1+x)|\leq x^2$ for any $x\in \ol{K}$ with $0<|x|\leq 1$. Since $\ol{K}\subseteq \R$, we can compute that $\ol{e}'(0)=1$. Now \cite[Theorem~3.3]{krappthesis} implies $(\ol{K},\ol{e})\subseteq (\R,\exp)$. Since $\exp$ is satisfies (GA) and (T$_{{1}}$), so does its restriction $\ol{e}$ to $\ol{K}$.
	\end{proof}

	Let $\mathcal{K}=(K,e,\ldots)$ be an o-minimal expansion of a real closed ordered (GAT$_1$)-exponential field $(K,e)$. Then $\mathcal{K}$ is called \textbf{exponentially bounded} if for any unary function $f\colon K\to K$ that is definable in $\mathcal{K}$ there are $n\in \N$ and $B\in K$ such that for any $x\in K$ with $x>B$ we have $|f(x)|<e^n(x)$. In this terminology, the question from the introduction can be phrased as follows:
	
	\begin{question}\thlabel{qu:expbdd}
		Is every o-minimal expansion of a real closed ordered (GAT$_1$)-exponential field exponentially bounded?
	\end{question}

	Let $(K,e)$ be a non-archimedean ordered exponential field with $v$-com\-pa\-ti\-ble exponential $e$. 
	Following the notation and terminology of \cite[Section~1.4]{kuhlmann}, we obtain the following maps: 
	The map $(-v\circ e|_{\bA})\colon \bA\to G, a\mapsto -v(e(a))$ is an order-preserving isomorphism of groups. Its inverse $(-v\circ e|_{\bA})^{-1}\colon G\to \bA$ is denoted by $h_e$. 
	Now $h_e$ induces an order-preserving bijection
	$$\widetilde{h}_e\colon \Gamma \to G^{<0}$$
	given by $\widetilde{h}_e(v_G(g))=v(h_e(g))$ for any $g\in G\setminus \{0\}$. The map $\widetilde{h}_e$ is called the \textbf{group exponential} on $G$ induced by $e$. The \textbf{natural contraction} $\chi_e$ (cf.\ \cite[Section~2.7]{kuhlmann}) induced by $e$ is given by
	$$\chi_e\colon G^{<0}\to G^{<0}, g\mapsto \widetilde{h}_e(v_G(g)).$$ 
	It has the property $\chi_e(v(a))=v(\ell(a))$ for any $a\in \bA^{>0}$. 
	Finally, the \textbf{shift map} $\zeta_e$ induced by $e$ is given by $$\zeta_e\colon \Gamma\to \Gamma, v_G(g)\mapsto v_G(\chi_e(g))$$ (where $g\in G^{<0}$).
	
	The natural contraction $\chi_e$ is one particular contraction map on $G^{<0}$. More generally, we call a monotonically increasing surjective map $\chi\colon G^{<0}\to G^{<0}$ a \textbf{contraction} if for any $g,h\in G^{<0}$ with $v_G(g)=v_G(h)$ we have $\chi(g)=\chi(h)$. Moreover, $\chi$ is called \textbf{centripetal} if for any $g\in G^{<0}$ it satisfies $g<\chi(g)$. Due to \cite[Corollary~2.29]{kuhlmann}, if $(K,e)$ is a non-archimedean ordered (GAT$_1$)-exponential field, then the natural contraction $\chi_e$ on $G^{<0}$ is centripetal.
	
	\begin{remark}
		Often in the literature, a contraction map $\chi$ is extended to $G$ symmetric about the origin, i.e.\ one sets $\chi(0)=0$ and $\chi(g)=-\chi(-g)$ for any $g\in G^{>0}$. Especially when the structure $(G,\chi)$ is studied model theoretically, $G$ has to be taken as the domain of $\chi$ (cf.\ \cite[Appendix~A]{kuhlmann}). Since we are interested in contractions induced by exponentials (and later transexponentials), in our context it suffices to take $G^{<0}$ as their domain. \qed
	\end{remark}
	
	We present further results on ordered exponential fields that will be used in later sections.
	
	\begin{lemma}\thlabel{lem:evgrowth}
		Let $(K,e)$ be a non-archimedean ordered (GAT$_1$)-exponential field.
		Let $x,y\in \PP_K$ with $v(x)<v(y)$. Then $v(e(x))<v(e(y))$.
	\end{lemma}

	\begin{proof}
		We have $x>2y$. Thus, $v(e(x))\leq v(e(2y))=2v(e(y))<v(e(y))$.
	\end{proof}

	The following is a direct consequence of \cite[Theorem~1.44]{kuhlmann} as well as \cite[Theorem~2.23]{kuhlmann} combined with \cite[Proposition~3.39]{krappthesis} and \Autoref{lem:t1exp}.

	\begin{lemma}\thlabel{fact:ccc}
		Suppose that $K$ is countable, non-archimedean and that $(K^{>0},\allowbreak \cdot,1,<)$ is divisible. Suppose further that $\ol{e}$ is a (GAT$_1$)-exponential on $\ol{K}$. Then the following are equivalent:
		\begin{enumerate}[label = (\roman*)]
			\item There is a (GAT$_1$)-exponential $e$ on $K$ with residue exponential $\ol{e}$.
			
			\item $G$ is isomorphic to $\coprod_{\Q}\ol{K}$ as an ordered group.
		\end{enumerate}
	\end{lemma}
		
\section{Theory of ordered transexponential fields}\label{sec:theory}
	
	\begin{definition}\thlabel{def:totf}
		We denote by $\Lotf$ the language $\{+,\cdot,0,1,<,e,T\}$ of or\-der\-ed transexponential fields, where $e$ and $T$ are unary function symbols. A structure $(K,+,\cdot,0,1,<,e,T)$ (from now on only abbreviated as $(K,e,T)$) is called an \textbf{ordered transexponential field} if it satisfies the following axioms:
	\begin{enumerate}[label = (\roman*)]
		\item $(K,+,\cdot,0,1,<)$ is an ordered field;
		
		\item $e$ is a (GAT$_1$)-exponential on $K$;
		
		\item $T\colon K\to K^{>0}$ is order-preserving and bijective;
		
		\item $T|_{[0,1]_K}=e|_{[0,1]_K}$;
		
		\item for any $x\in K^{>0}$ we have $T(x+1)=e(T(x))$;
		
		\item for any $x\in K^{<0}$ we have $T(x)T(-x)=1$.
	\end{enumerate}
	Moreover, we say that $T$ is a \textbf{transexponential} on $(K,e)$ or that $T$ is a transexponential on $K$ \textbf{compatible} with $e$. 
	Its inverse $T^{-1}$ is called a \textbf{cislogarithm} and we denote it by $L$. 
	The $\Lotf$-theory of ordered transexponential fields axiomatised by the above is denoted by $T_{\mathrm{otf}}$.
	\end{definition}
	
	For any (GAT$_1$)-exponential $e$, the growth axiom scheme ensures that $e$ eventually exceeds any polynomial.
	However, in general this growth property fails in the absence of (GA) (see, e.g., \cite[Section~2.6]{kuhlmann} for a more detailed discussion).
	We will see in \Autoref{ex:nogrowth} that the axiomatisation of ordered transexponential fields does likewise not ensure that $T$ eventually exceeds every iterate of $e$.
	
	\begin{definition}\thlabel{def:grwoth}
		Let $(K,e,T)$ be an ordered transexponential field. Then $T$ has the \textbf{growth property} if for any $n\in \N$ and for any $a\in K$ with $a\geq (n+1)^2$ we have $T(a)>e^n(a)$.
	\end{definition}

	\begin{example}\thlabel{ex:realtrans}
		Recall that we denote by $(\R,\exp)$ the real exponential field, where $\exp$ denotes the standard exponential $x\mapsto \mathrm{e}^x$. 
		Note that due to \Autoref{lem:t1exp}, $\exp$ is the only (GAT$_1$)-exponential on $\R$. 
		The only transexponential on $\R$ compatible with $\exp$ is given by
		$$T_{\R}(x)=\exp^{m+1}(x-m)$$
		for $x\in \R^{\geq0}$ and $m\in \N\cup\{0\}$ with $x\in [m,m+1)_{\R}$, as well as $T_{\R}(x)=(T_{\R}(-x))^{-1}$ for $x\in \R^{<0}$. We will see in \Autoref{prop:archgrowth} that $T_{\R}$ has the growth property. \qed
	\end{example}

	\Autoref{ex:realtrans} presents our main example of an archimedean ordered transexponential fields. Specific examples of non-archimedean ordered transexponential fields will be constructed in \Autoref{sec:countable}.	
	We now establish some basic properties of transexponentials, which will be used to show in \Autoref{thm:notomin} that $T_{\mathrm{otf}}$ does not admit o-minimal models.
	
	\begin{lemma}\thlabel{lem:basics}
		Let $(K,e,T)$ be an ordered transexponential field. Then the following hold:
		\begin{enumerate}[label = (\roman*)]
			\item\label{lem:basics:1} For any $x\in K^{\geq 0}$ and $m\in \Z$ with $x+m\geq 0$ we have $T(x+m)=e^m(T(x))$.
			
			\item\label{lem:basics:2} For any $x\in K^{>1}$, if $T$ is differentiable in $x$, then $$T'(x) = T(x)T'(x-1) = e(T(x-1))T'(x-1).$$
			
			\item\label{lem:basics:3} For any $m\in \N\cup \{0\}$ and any $x\in [m,m+1]_K$ we have $$T'(x)=e(x-m)\cdot \prod_{i=1}^me(T(x-i)).$$
			In particular, $T$ is continuously differentiable on $R_v$. 
		\end{enumerate}
	\end{lemma}

	\begin{proof}
			\begin{enumerate}[label = (\roman*)]
			\item If $m\geq 0$, then we can iteratively apply that $T(y+1)=e(T(y))$ holds for any $y\in K^{\geq 0}$. If $m<0$, then $T(x)=T(x+m-m)=e^{-m}(T(x+m))$ and thus
			$T(x+m)=\ell^{-m}(T(x))=e^m(T(x))$.
			
			\item By \Autoref{lem:t1exp}, we have $e'=e$. The claim follows by the chain rule, which holds for differentiable functions on any ordered field.
			
			\item 	For any $x\in [m,m+1]_K$,  \autoref{lem:basics:1} implies $T(x)=e^m(T(x-m))=e^{m+1}(x-m)$.
			First note that for any $x\in [0,1)_K$ we have $T'(x)=e(x)$.
			Now let $m\geq 1$.
			By iterated application of the chain rule, we obtain for $x\in (m,m+1)_K$ that
			\begin{align*}
				T'(x)&=e^{m+1}(x-m)\cdot e^m(x-m)\cdot\ldots\cdot e^2(x-m)\cdot e(x-m)\\
				&= e(x-m)\cdot e^{m}(T(x-m))\cdot e^{m-1}(T(x-m))\cdot\ldots\cdot e(T(x-m))\\
				&= e(x-m)\cdot e(T(x-1))\cdot e(T(x-2))\cdot\ldots\cdot e(T(x-m))\\
				&= e(x-m)\cdot \prod_{i=1}^me(T(x-i)).
			\end{align*}
			To prove the claim in the interval boundaries of $[m,m+1]_K$, we simply note that $T'(x)$ can be continuously extended to any $n\in \N$. Indeed, taking left- and right-sided limits, we obtain
			\begin{align*}
				\lim\limits_{x\searrow n} T'(x)&=e(n-n)\cdot \prod_{i=1}^ne(T(n-i))\\&=e(T(0))\cdot \prod_{i=1}^{n-1}e(T(n-i))=\lim\limits_{x\nearrow n} T'(x).
			\end{align*}
			We have thus established that $T$ is continuously differentiable on $R_v^{\geq 0}$. To obtain that $T$ is continuosly differentiable on $R_v$, simply note that for $x\in K^{\leq 0}$ we have
			$$T'(x)=\frac{T'(-x)}{(T(-x))^2}.$$
		\end{enumerate}
	\end{proof}

	\begin{remark}
		The properties of an exponential already ensure that it is differentiable everywhere if and only if it is differentiable in at least one point (cf.\ \cite[Corollary~2.14]{krappthesis}). This is different for transexponentials: 
		while \Autoref{lem:basics}~\autoref{lem:basics:2} expresses the derivative of a transexponential in the points where it is differentiable and \Autoref{lem:basics}~\autoref{lem:basics:3} establishes that any transexponential is differentiable on $R_v$, there are examples of transexponentials that are not differentiable in some points outside $R_v$ (see \Autoref{rmk:notgrow}).
		%Use in example a map 0 to a, 1/2 to a+1, 1 to exp(a). this is not differentiable in 1/2
		\qed
	\end{remark}

	\begin{theorem}\thlabel{thm:notomin}
		Any $(K,e,T)\models T_{\mathrm{otf}}$ is not o-minimal. 
	\end{theorem}

	\begin{proof}
		Let $A\subseteq K^{>0}$ be the $\Lotf$-definable set consisting of all $x\in K^{>0}$ in which $T$ is not twice continuously differentiable. We show that $A\cap R_v = \N$, yielding that $A$ cannot be expressed as a finite union of points and open intervals in $K$.
				
		Let $m\in \N\cup \{0\}$ and let $x\in (m,m+1)_K$.	By \Autoref{lem:basics}~\autoref{lem:basics:3}, we have $$T'(x)=e(x-m)\cdot \prod_{i=1}^me(T(x-i)).$$
		Thus, applying the product rule we obtain,
		$$T''(x)=T'(x)\brackets{1+\sum_{i=1}^mT'(x-i)}.$$
		Since $T'$ is continuously differentiable on $(m,m+1)_K$, on this open interval $T$ is twice continuously differentiable.
		Now let $n\in \N$ and consider left- and right-sided limits of $T''$ as $x$ approaches $n$. We obtain
		\begin{align*}\lim\limits_{x\searrow n} T''(x)&=T'(n)\brackets{1+\sum_{i=1}^nT'(n-i)} \text{ and}\\
		\lim\limits_{x\nearrow n} T''(x)&=T'(n)\brackets{1+\sum_{i=1}^{n-1}T'(n-i)}.
		\end{align*}
		Hence, $\lim\limits_{x\searrow n} T''(x)\neq \lim\limits_{x\nearrow n} T''(x)$, showing that $T$ is not twice continuously differentiable in $n$.
	\end{proof}

	\begin{corollary}\thlabel{cor:notomin}
		Let $\mathcal{K}=(K,e,\ldots)$ be an o-minimal expansion of a real closed ordered (GAT$_1$)-exponential field $(K,e)$. Then for any definable function $f\colon K\to K$ in $\mathcal{K}$ we have $(K,e,f)\not\models T_{\mathrm{otf}}$. 
	\end{corollary}

	\begin{proof}
		Let $f\colon K\to K$ be definable in $\mathcal{K}$. Then $(K,e,f)$ is o-minimal, and \Autoref{thm:notomin} implies $(K,e,f)\not\models T_{\mathrm{otf}}$.
	\end{proof}

	\begin{example}\thlabel{ex:nondefinable}
		\Autoref{cor:notomin} implies that $T_{\R}$ is not definable in the o-minimal structure $(\R,\exp)$. 
		More generally, for any o-minimal ordered exponential field $(K,e)$ whose exponential satisfies $e'=e$, \cite[Proposition~4.36]{krappthesis} shows that $e$ is a (GAT$_1$)-exponential. Thus, \Autoref{cor:notomin} yields that $(K,e)$ admits no definable transexponential compatible with $e$.
		\qed
	\end{example}

	\begin{remark}
		\Autoref{ex:nondefinable} stands in contrast to the following observation: 
		Given an ordered transexponential field $(K,e,T)$, the exponential $e$ is already definable in the structure $(K,T)$. Indeed, let $x\in K$. If $x>1$, then $$e(x)=T(L(x)+1),$$ if $0\leq x\leq 1$, then $e(x)=T(x)$, and if $x<0$, then $e(x)=(e(-x))^{-1}=(T(L(-x)+1))^{-1}$.
		\qed
	\end{remark}

	One approach towards answering \Autoref{qu:expbdd} affirmatively would now be as follows: 
	Let $\mathcal{K}$ be an o-minimal expansion of a real closed ordered (GAT$_1$)-exponential field $(K,e)$. Assume, for a contradiction, that $\mathcal{K}$ is not exponentially bounded. 
	If one can now define in $\mathcal{K}$ a function $T\colon K\to K$ such that $(K,e,T)\models  T_{\mathrm{otf}}$, this leads to the required contradiction by \Autoref{cor:notomin}.
	
\section{Induced structure on the residue field}\label{sec:res}

	In order to begin our systematic study of the behaviour of transexponentials on the valuation theoretic invariants of $K$, we first consider $T$ restricted to the valuation ring $R_v$. 
	\textbf{Throughout this section, we let $(K,e,T)$ be an ordered transexponential field.}	
	
		\begin{proposition}\thlabel{prop:transrv}
		\begin{enumerate}[label = (\roman*)]
			\item $T|_{I_v}$ is an order-preserving isomorphism from $(I_v,\allowbreak +,0,<)$ to $(1+I_v,\cdot,1,<)$.
			
			\item\label{prop:transrv:2} $T|_{\OO_v}$ is an order-preserving bijection from $(\OO_v,<)$ to $(\U_v^{>0},<)$.
		\end{enumerate}
	\end{proposition}
	
	\begin{proof}		
		\begin{enumerate}[label = (\roman*)]
			\item Simply note that $T|_{I_v}=e|_{I_v}$, which has the required property by $v$-compatibility. 
			
			\item We only need to show that $T(\OO_v)=\U_v^{>0}$. First we have
			\begin{align*}
			T(\OO_v^{\geq 0})&=T\!\brackets{\bigcup_{n\in \N}[n-1,n)_K}=\bigcup_{n=1}^\infty T([n-1,n)_K)=\bigcup_{n=1}^\infty [e^n(0),e^n(1))_K\\
			&=[e(0),\infty)_{R_v} = [1,\infty)_{R_v}=\mathcal{U}_v^{\geq 1}.
			\end{align*}
			Hence,
			\begin{align*}
			T(R_v)&=T(R_v^{\geq 0})\cup T(R_v^{\leq  0})\\&=\{x\mid v(x)=0\text{ and }x\in K^{\geq1}\}\cup \{x^{-1}\mid v(x)=0\text{ and }x\in K^{\geq1}\}\\&=\mathcal{U}_v^{>0}.
			\end{align*}                                      
		\end{enumerate}
	\end{proof}

	We now show in \Autoref{prop:rvgrowth} below that on $R_v$, a transexponential eventually exceeds any iterate of $e$. 
	For any $a\in R_v$, we denote by $\lfloor a \rfloor \in \Z$ the integer part of $a$, i.e.\ the unique integer $z\in \Z$ with $z\leq a<z+1$.
	
	\begin{lemma}\thlabel{lem:growth1}
		For any $a\in \OO_v^{\geq 0}$, we have $T(a)=e^{\lfloor a\rfloor+1}(a-\lfloor a \rfloor)$. In particular, if $a\in \N$, then $T(a)=e^{a+1}(0)=e^a(1)$. 
	\end{lemma}

	\begin{proof}
		We simply apply \Autoref{lem:basics}~\autoref{lem:basics:1} to $x=a-\ip{a}$ and $m=\ip{a}$. 
	\end{proof}

	Note that \Autoref{lem:growth1} implies that for any transexponential on $K$ compatible with $e$, the values it takes on $R_v$ are uniquely determined. 

	\begin{lemma}\thlabel{lem:tgrowthov1}
		For any $n\in \N$, we have $T(n)>2^n$.
	\end{lemma}

	\begin{proof}
		Since $e(4)> 4^2=16$, we obtain $$e(1)=\sqrt[4]{e(4)}>\sqrt[4]{16}=2.$$
		To show that $T(n)>2^n$, first note that $T(0)=e(0)=1=2^0$ and then argue inductively by using
		$$T(n)=e(T(n-1))\geq e(2^{n-1})>2^{2^{n-1}}\geq 2^n.$$
	\end{proof}

	\begin{lemma}\thlabel{lem:tgrowthov2}
		For any $n,m\in \N$ with $m\geq (n+1)^2$ we have $2^{m-n}\geq m + 1$.
	\end{lemma}

	\begin{proof}
		Since $m-n-2\geq 1$, we obtain
		\begin{align*}
			2^{m-n}&=4\cdot 2^{m-n-2}\geq 4(1+m-n-2)\\
			&\geq (m-n-1)+3((n+1)^2-n-1)
			=m+3n^2+2n+2>m+1.
		\end{align*}
	\end{proof}
	
	\begin{proposition}\thlabel{prop:rvgrowth}
		For any $n\in \N$ and any $x\in \OO_v$ with $x\geq (n+1)^2$ we have $T(x)> e^n(x)$.
	\end{proposition}

	\begin{proof}
		Let $m=\lfloor x\rfloor$.
		By \Autoref{lem:growth1}, \Autoref{lem:tgrowthov1} and \Autoref{lem:tgrowthov2}, we obtain
		\begin{align*}
			T(x)&\geq T(m) = e^n(T(m-n)) > e^n(2^{m-n})\\
			&\geq e^n(m+1)\geq e^n(x).
		\end{align*}
	\end{proof}

	Just as $v$-compatibility of $e$ ensures that $e$ induces a residue exponential $\ol{e}$ on $\ol{K}$, we now show in \Autoref{thm:restrans} below that also $T$ induces a residue transexponential $\ol{T}\colon \ol{K}\to \ol{K}^{>0}$.

	\begin{lemma}\thlabel{lem:valprop}
		For any $a,b\in \OO_v$ with $a-b\in \II_v$, we have $T(a)-T(b)\in \II_v$.
	\end{lemma}

	\begin{proof}
		We may assume that $a<b$. 
		If $a<0<b$, then $a,b\in I_v$. Thus, $T(a)=e(a)$, $T(b)=e(b)$ and $T(a),T(b)\in 1+I_v$, yielding the claim.	
		
		Now suppose that $0<a<b$. 
		Let $m=\ip{a}$. If also $m=\ip{b}$, then 
		\begin{align*}
			T(a)-T(b)&=e^{m+1}(a-m)-e^{m+1}(b-m).
		\end{align*}
		Now $(a-m)-(b-m)\in \II_v$. We can iteratively apply the $v$-compatiblity of $e$ to obtain the desired conclusion.
		Thus, suppose that $\ip{b}=m+1$. Then $a<m+1\leq b$ and
		$$T(a)-T(b)=e^{m+1}(a-m)-e^{m+1}(e(b-m-1)).$$
		Now $a-m\in 1+I_v$ and $b-m-1\in I_v$, so also $e(b-m-1)\in 1+I_v$. Hence, $(a-m)-e(b-m-1)\in I_v$ and, again, we can apply that $e$ is $v$-compatible.
		
		Finally, suppose that $a<b<0$. As $v(T(-a)T(-b))\leq 0$, we obtain
		\begin{align*}
			v(T(a)-T(b))&=	v(T(-a)^{-1}-T(-b)^{-1}) = v\!\brackets{\frac{T(-b)-T(-a)}{T(-a)T(-b)}}\\
			&\geq v(T(-b)-T(-a))>0.
		\end{align*}
	\end{proof}

	\begin{theorem}\thlabel{thm:restrans}
		The map $\ol{T}\colon \ol{K}\to \ol{K}^{>0}$ given by $\ol{a}\mapsto \ol{T(a)}$ for any $a\in R_v$ is a transexponential on $\ol{K}$ which is compatible with $\ol{e}$.
	\end{theorem}

	\begin{proof}
		By \Autoref{lem:gatgat}, $\ol{e}$ is a (GAT$_1$)-exponential on $\ol{K}$.
		\Autoref{lem:valprop} implies that the map $\ol{T}$ is well-defined. By \Autoref{prop:transrv}~\autoref{prop:transrv:2}, it is an order-preserving bijection from $\{\ol{a}\mid a\in R_v\}=\ol{K}$ to $\{\ol{a}\mid a\in \mathcal{U}_v^{>0}\}=\ol{K}^{>0}$. 
		Since $T|_{[-1,1]_K}=e|_{[-1,1]_K}$, we obtain $$\ol{T}(\ol{a})=\ol{T(a)}=\ol{e(a)}=\ol{e}(\ol{a})$$ for any $a\in R_v$ with $\ol{a}\in [0,1]_{\ol{K}}$. Finally, let $a\in R_v$ with $\ol{a}>0$. Then
		$\ol{T}(\ol{a}+1)=\ol{T(a+1)}=\ol{e(T(a))}=\ol{e}(\ol{T}(\ol{a}))$, and
		$\ol{T}(-\ol{a})\ol{T}(\ol{a})=\ol{T(-a)T(a)}=1$, as required.
	\end{proof}

	\Autoref{thm:restrans} gives rise to a class of archimedean ordered transexponential fields. Namely, if $(K,e,T)$ is a non-archimedean model of $T_{\mathrm{otf}}$, then $(\ol{K},\ol{e},\ol{T})$ becomes an archimedean model of $T_{\mathrm{otf}}$. We show in \Autoref{thm:arch} below that any archimedean model of $T_{\mathrm{otf}}$ is a substructure of $(\R,\exp,T_{\R})$ introduced in \Autoref{ex:realtrans}.
	Note that for any archimedean ordered field $K$ we have $R_v=K=\ol{K}$.

	\begin{proposition}\thlabel{prop:archgrowth}
		Let $(K,e)$ be an archimedean ordered (GAT$_1$)-exponential field. Then there is a unique transexponential $T$ on $K$ that is compatible with $e$. Moreover, $T$ has the growth property.
	\end{proposition}

	\begin{proof}
		\Autoref{lem:growth1} implies that the only transexponential $T$ on $K$ compatible with $e$ is given by
		$$T\colon K\to K^{>0}, x\mapsto \begin{cases}
		e^{\lfloor x\rfloor+1}(x-\lfloor x \rfloor)&\text{ if }x\geq 0,\\
		(e^{\lfloor x\rfloor+1}(x-\lfloor x \rfloor))^{-1}&\text{ if }x< 0.
		\end{cases}$$
		It follows immediately from \Autoref{prop:rvgrowth} that $T$ has the growth property.
	\end{proof}
	
	\begin{remark}
		\Autoref{prop:archgrowth} shows that the class of archimedean ordered fields that can be endowed with a (GAT$_1$)-exponential coincides with the class of archimedean ordered fields that can be expanded to an ordered transexponential field. 
		Hence, an archimedean ordered field $K$ can be expanded to a model of $T_{\mathrm{otf}}$ if and only if it is exponentially closed in $(\R,\exp)$. Necessary and sufficient conditions under which countable non-archimedean ordered fields can be expanded to a model of $T_{\mathrm{otf}}$ will be presented in \Autoref{thm:nonarchexpand}. \qed
	\end{remark}

	\begin{proposition}\thlabel{thm:arch}
		Let $(K,e,T)$ be an archimedean ordered transexponential field. Then $(K,e,T)\subseteq (\R,\exp,T_{\R})$. More precisely, $(K,e,T)$ embeds into $(\R,\exp,T_{\R})$ as an $\Lotf$-structure via a unique embedding.
	\end{proposition}

	\begin{proof}
		By \Autoref{lem:gatgat}, $(K,e)\subseteq (\R,\exp)$. It remains to note that the transexponential $T$ from \Autoref{prop:archgrowth} coincides with $T_{\R}$ restricted to $K$.
	\end{proof}

	As a consequence of \Autoref{prop:archgrowth} and \Autoref{thm:arch}, we obtain the following corollary to \Autoref{thm:restrans}.

	\begin{corollary}
		Let $(K,e,T)$ be a non-archimedean ordered transexponential field. Then the transexponential $\ol{T}$ on $(\ol{K},\ol{e})$ has the growth property. Moreover, $(\ol{K},\ol{e},\ol{T})\subseteq (\R,\exp,T_{\R})$.
	\end{corollary}
	
\section{Left-transexponential on the positive infinite part}\label{sec:pos}

	Let $(K,e,T)$ be a non-archimedean ordered transexponential field.
	The behaviour of $T$ on the set of finite elements $R_v$ was described in 
	\Autoref{prop:transrv}~\autoref{prop:transrv:2}; namely, $T$ maps $R_v$ bijectively onto $\mathcal{U}_v^{>0}$. Since $T$ is an order-preserving bijection from $K$ to $K^{>0}$, this implies that $T$ maps the set of negative infinite elements $-\mathbf{P}_K=\{x\in K^{<0}\mid v(x)<0\}$ onto the set of positive infinitesimals $I_v^{>0}$ and the set of positive elements $\mathbf{P}_K=\{x\in K^{>0}\mid v(x)<0\}$ onto itself. %We obtain the following picture: PICTURE
	%\texttt{Include picture of ordered line of $K$ here.}
	We now focus on the restriction of $T$ to $\mathbf{P}_K$. Inspired by the decomposition of exponentials into left-, middle- and right-exponentials from \cite[page~24]{kuhlmann}, where the left-exponential is defined on the complement to the valuation ring, we introduce the following terminology.
	\begin{definition}
		Let $(K,e,T)$ be a non-archimedean ordered transexponential field. Then the map
		$\TL:=T|_{\mathbf{P}_K}$ is called a \textbf{left-transexponential} on $(K,e)$.
	\end{definition}
	Note that any left-transexponential is an order-preserving bijection from $\mathbf{P}_K$ to $\mathbf{P}_K$. 
	%Since $T(\OO_v)=\U_v^{>0}$, we obtain that for any $a\in K^{>0}$ with $v(a)<0$ also $v(T(a))<0$. 
	
	\begin{remark}\thlabel{rmk:alltransexp}
		Let $(K,e)$ be a non-archimedean ordered (GAT$_1$)-exponential field. Due to \Autoref{lem:growth1}, for any transexponential $T$ on $K$ compatible with $e$, the restriction of $T$ to $R_v$ is uniquely determined. Since $K$ decomposes as an ordered set into $-\PP_K$, $R_v$ and $\PP_K$, and the values $T$ takes on $-\PP_K$ are determined by the values it takes on $\PP_K$, the only freedom of choice we have to construct $T$ is its behaviour on $\PP_K$. More precisely, given a transexponential $T$ on $(K,e)$ with corresponding left-transexponential $\TL=T|_{\mathbf{P}_K}$, we have
		\begin{align}
			T\colon K\to K^{>0}, x\mapsto \begin{cases}
		(\TL(-x))^{-1}	&\text{ if } x\in -\PP_K,\\
		(e^{\lfloor x\rfloor+1}(x-\lfloor x \rfloor))^{-1}	&\text{ if } x\in R_v^{<0},\\
		e^{\lfloor x\rfloor+1}(x-\lfloor x \rfloor)	&\text{ if } x\in R_v^{\geq0},\\
		\TL(x)	&\text{ if } x\in \PP_K.
		\end{cases}
		\label{eq:transform}
		\end{align} 
	\qed
	\end{remark}

	From the discussion above, we immediately obtain the following.

	\begin{lemma}\thlabel{lem:lefttransax}
		Let $(K,e)$ be a non-archimedean ordered (GAT$_1$)-exponential field. Let $\TL\colon \PP_K\to \PP_K$ be an arbitrary map. Then the map $T$ given by (\autoref{eq:transform}) is a transexponential on $(K,e)$ if and only if the following hold:
		\begin{enumerate}[label = (\roman*)]
			\item\label{lem:lefttransax:1} $\TL$ is an order-preserving bijection from $\PP_K$ onto $\PP_K$.
			
			\item\label{lem:lefttransax:2} For any $x\in \PP_K$, we have $\TL(x+1)=e(\TL(x))$.
		\end{enumerate}
	\end{lemma}
	
	Given a non-archimedean field $K$ and a (GAT$_1$)-exponential $e$ on $K$, recall that $K$ additively decomposes as $K = \bA\amalg R_v$ and that for $a,b\in K$ we write $a\sim_e b$ if there exists $n\in \N$ with $a<e^n(b)$ and $b<e^n(a)$.
	By \cite[Corollaries~3.2~and~3.34]{kuhlmann}, the order-type of the ordered set $\PP_K/\!\sim_e$ is the \textbf{principal exponential rank} of $(K,e)$.	By abuse of terminology, we also call the ordered set $\PP_K/\!\sim_e$ the principal exponential rank of $(K,e)$.
	
	In the following, we present a construction method (\Autoref{cons:trans}) for transexponentials on $(K,e)$, provided that $\bA^{>0}$ is isomorphic to $\PP_K/\!\sim_e$ as an ordered set. Subsequently, we prove that any transexponential on $(K,e)$ can be obtained via this method (\Autoref{thm:alltransexp}).

	\begin{construction}\thlabel{cons:trans}
		Let $K$ be a non-archimedean ordered field and let $e$ be a (GAT$_1$)-exponential on $K$. Suppose that there exists an order-preserving bijection
		$$\varphi\colon \bA^{>0}\to \PP_K/\!\sim_e.$$
		For any $a\in \bA^{>0}$, let $c_a\in \varphi(a)$ and let $f_a\colon [0,1]_K\to [c_a,e(c_a)]_K$ be an order-preserving bijection. For instance, one can choose  $f_a\colon t\mapsto (1-t)c_a+te(c_a)$.
		We now construct a transexponential $T$ on $(K,e)$ depending on the family $(c_a,f_a)_{a\in \bA^{>0}}$. 
		
		Noting that $$\PP_K=\bigcupdot_{a\in \bA^{>0}}(a+R_v),$$
		(see the picture below),
		we define an order-preserving bijection $$\TL\colon \PP_K\to \PP_K$$
		as follows: for any $a\in \bA^{>0}$ and any $b\in R_v$, we set 
		$$\TL(a+b)=e^{\ip{b}}(f_a(b-\ip{b})).$$
		
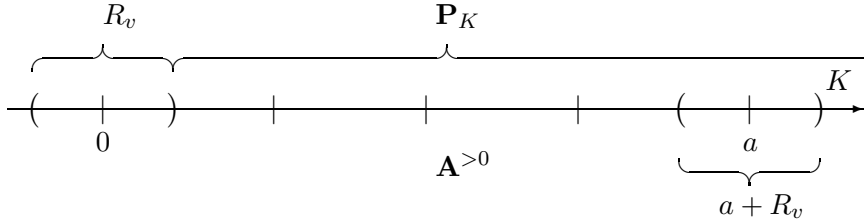
					\begin{center}
			\setlength{\unitlength}{0.002\textwidth}
			\begin{picture}(500,140)(-70,-70)
			
			%line
			\put(-50,0){\vector(1,0){450}}
			\put(380,10){{$K$}}
			
			%zero
			\put(-2,-4){{$|$}}
			\put(-3.5,-22){{$0$}}
			
			%valuation ring
			\put(-40,-4){{$\left(\hspace{24pt}\rule{0pt}{10pt}
					\hspace{24pt}\right)$}}
			\put(-30,20){\oval(14,14)[tl]}
			\put(30,20){\oval(14,14)[tr]}
			\put(-7,34){\oval(14,14)[br]}
			\put(7,34){\oval(14,14)[bl]}
			\put(-30,27){\line(1,0){23}}
			\put(7,27){\line(1,0){23}}	
			\put(0,45){{$R_v$}}
			
			%infinite elements
			\put(44,20){\oval(14,14)[tl]}
			\put(44,27){\line(1,0){130}}
			\put(174,34){\oval(14,14)[br]}
			\put(188,34){\oval(14,14)[bl]}
			\put(188,27){\line(1,0){212}}
			\put(175,45){{{\bf P}$_K$}}
			
			%A
			\put(175,-35){$\bA^{>0}$}
			
			%a
			\put(338,-4){{$|$}}
			\put(336.5,-22){{$a$}}
			\put(300,-4){{$\left(\hspace{24pt}\rule{0pt}{10pt}
				\hspace{24pt}\right)$}}
			
			%a + R_v
			\put(310,-26){\oval(14,14)[bl]}
			\put(370,-26){\oval(14,14)[br]}
			\put(333,-40){\oval(14,14)[tr]}
			\put(347,-40){\oval(14,14)[tl]}
			\put(310,-33){\line(1,0){23}}
			\put(347,-33){\line(1,0){23}}	
			\put(324,-55){{$a+R_v$}}
			
			%a's
			\put(248,-4){{$|$}}
			\put(168,-4){{$|$}}
			\put(88,-4){{$|$}}
				
			%\put(52,10){{$\ldots$}}
			%\put(100,0){{$\left(\hspace{20pt}\rule{0pt}{10pt}
			%		\hspace{20pt}\right)$}}
			%\put(152,10){{$\ldots$}}
			%\put(250,0){{$\left(\hspace{52pt}\rule{0pt}{10pt}
			%		\hspace{52pt}\right)$}}
			%\put(347,10){{$\ldots$}}
			%\put(220,15){{$a^{1/n}$}}
			%\put(250,0){{$|$}}
			%\put(250,15){{$a$}}
			%\put(280,0){{$|$}}
			%\put(280,15){{$a^n$}}
			%\put(190,10){{$\ldots$}}
			%\put(310,10){{$\ldots$}}

			%\put(180,27){\line(1,0){63}}
			%\put(320,27){\line(-1,0){63}}
			%\put(180,20){\oval(14,14)[tl]}
			%\put(320,20){\oval(14,14)[tr]}
			%\put(243,34){\oval(14,14)[br]}
			%\put(257,34){\oval(14,14)[bl]}
			%\put(250,45){{$[a]^.$}}

			\end{picture}
		\end{center}
		
		The transexponential $T$ on $(K,e)$ with left-transexponential $\TL$ is now given by the formula  (\autoref{eq:transform}). 
		In order to show that $T$ is indeed a transexponential on $K$ compatible with $e$, by \Autoref{lem:lefttransax} it remains to verify the following.
		\begin{enumerate}[label = (\roman*)]
			\item \emph{$\TL$ is an order-preserving bijection from $\PP_K$ onto $\PP_K$:}
				Let $x,y\in \PP_K$ with $x<y$, let $a,a'\in \bA^{>0}$ and $b,b'\in R_v$ such that $x=a+b$ and $y=a'+b'$. Recall that $e|_{\PP_K}$ is an order-preserving bijection from $\PP_K$ to $\PP_K$. Since $\TL(x)=e^{\ip{b}}(f_a(b-\ip{b}))$, we obtain $\TL(x)\in \PP_K$. Indeed, it suffices to note that $f_a(b-\ip{b})\in[c_a,e(c_a)]_K\subseteq \PP_K$, as $c_a\in \PP_K$. 
				
				In order to show that $\TL$ preserves the order on $\PP_K$, first note that $\TL(x)\in [c_a]_e=\varphi(a)$ and $\TL(y)\in [c_{a'}]_e=\varphi(a')$. If $a<a'$, then $\varphi(a)<\varphi(a')$, yielding that $\TL(x)<\TL(y)$. If $a=a'$, then $b<b'$. If $\ip{b}=\ip{b'}$, then
				$$\TL(x)=e^{\ip{b}}(f_a(b-\ip{b}))<e^{\ip{b}}(f_a(b'-\ip{b}))=\TL(y),$$
				as $f_a$ is order-preserving.
				Otherwise, we have
				$$f_a(b-\ip{b})<e(c_a)\leq e^{\ip{b'}-\ip{b}}(c_a)\leq e^{\ip{b'}-\ip{b}}(f_a(b'-\ip{b'})).$$
				Thus,
				$$\TL(x)=e^{\ip{b}}(f_a(b-\ip{b}))<e^{\ip{b'}}(f_a(b'-\ip{b'}))=\TL(y).$$

			\item \emph{For any $x\in \PP_K$, we have $\TL(x+1)=e(\TL(x))$:} Let $a\in \bA^{>0}$ and $b\in R_v$ such that $x=a+b$. Thus
			$$\TL(x+1)=e^{\ip{b+1}}(f_a(b+1-\ip{b+1}))=e^{\ip{b}+1}(f_a(b-\ip{b}))=e(\TL(x)).$$

		\end{enumerate}		
		\qed
	\end{construction}
		
	The following complements \Autoref{prop:archgrowth}, which showed that in the archimedean case a transexponential always exists.
	
	\begin{proposition}\thlabel{prop:condextransexp}
		Let $(K,e)$ be a non-archimedean ordered (GAT$_1$)-ex\-po\-nen\-tial field. Then $(K,e)$ admits a transexponential if and only if $\mathbf{A}^{>0}$ is isomorphic as an ordered set to $\PP_K/\!\sim_e$.
	\end{proposition}

	\begin{proof}
		If $\mathbf{A}^{>0}$ is isomorphic to $\PP_K/\!\sim_e$, then \Autoref{cons:trans} shows how a left-transexponential $\TL$ and subsequently a transexponential on $(K,e)$ can be constructed. Conversely, suppose that $T$ is a transexponential on $(K,e)$ with left-transexponential $\TL$. We show that the map 
		$$\varphi\colon \bA^{>0}\to \PP_K/\!\sim_e, a\mapsto [\TL(a)]_e$$
		is an order-preserving bijection. First let $a,a'\in \bA^{>0}$ with $a<a'$. Then for any $n\in \N$ we have $a+n<a'$. Hence,
		$e^n(\TL(a))=\TL(a+n)<\TL(a')$, showing that $[\TL(a)]_e<[\TL(a')]_e$. In order to verify that $\varphi$ is a bijection, let $c\in \PP_K$ and set $d=\TL^{-1}(c)\in \PP_K$. Then there are $a\in \bA^{>0}$ and $b\in R_v$ such that $d=a+b$. Then
		$$\TL(d)\leq e^{\ip{b}+1}(\TL(a))\text{ and }\TL(a)\leq e^{\ip{-b}+1}(\TL(d)),$$
		yielding $\TL(a)\sim_e\TL(d)$. 
		Hence,
		$$\varphi(a)=[\TL(a)]_e=[\TL(d)]_e=[c]_e,$$
		as required.
	\end{proof}

	\begin{corollary}\thlabel{cor:dlo}
		Let $(K,e,T)$ be a non-archimedean ordered transexponential field. Then the principal exponential rank of $(K,e)$ is a dense linear order without endpoints.
	\end{corollary}

	\begin{proof}
		This follows from \Autoref{prop:condextransexp}, as $\bA^{>0}$ is a dense linear order without endpoints.
	\end{proof}
	
	We conclude this section by showing that any transexponential on a non-archimedean ordered exponential field can already be obtained by the method presented in \Autoref{cons:trans}.
	
	\begin{theorem}\thlabel{thm:alltransexp}
		Let $(K,e,T)$ be a non-archimedean ordered transexponential field. Then there exist
		an order-preserving bijection $\varphi\colon \bA^{>0}\to \PP_K/\!\sim_e$ and, for each $a\in \bA^{>0}$, some $c_a\in \varphi(a)$ and an order-preserving bijection $f_a\colon [0,1]_K\to [c_a,e(c_a)]_K$ such that $T$ is obtained by \Autoref{cons:trans} applied to the family $(c_a,f_a)_{a\in \bA^{>0}}$. 
	\end{theorem}

	\begin{proof}
		As in the proof of \Autoref{prop:condextransexp}, we set
		$$\varphi\colon \bA^{>0}\to \PP_K/\!\sim_e, a\mapsto [\TL(a)]_e.$$
		For any $a\in \bA$ we set 
		$$c_a=\TL(a)\text{ and }f_a\colon [0,1]_K\to [c_a,e(c_a)], t\mapsto \TL(a+t).$$
		Let $x\in \bA^{>0}$ and let $a\in \bA^{>0}$ and $b\in R_v$ with $x=a+b$. Then
		$$\TL(x)=\TL(a+b)=e^{\ip{b}}(\TL(a+b-\ip{b}))=e^{\ip{b}}(f_a(b-\ip{b})).$$
		Hence, $\TL$ and consequently $T$ are obtained by applying \Autoref{cons:trans} to $(c_a,f_a)_{a\in \bA^{>0}}$. 
	\end{proof}

\section{Induced structure on the value group}\label{sec:val}
	
	\textbf{Throughout this section, unless otherwise specified, let $(K,e,T)$ be a non-archimedean or\-dered transexponential field.} Recall that we denote by $G$ the value group of $K$ under the natural valuation $v$ and by $\Gamma$ the rank of $G$, i.e.\ the ordered value set under the natural valuation $v_G$ on $G$. 
	This section initiates the study of the map which a transexponential or, more precisely, the corresponding cislogarithm induces on $G$. 
	
	\begin{lemma}\thlabel{lem:valpres}
		Let $a,b\in \mathbf{P}_K$ with $b-a\geq 1$. Then $v(T(a))>v(T(b))$.
	\end{lemma}
	
	\begin{proof}
		Let $n\in \N$. Since $T(a)>n$ and $e$ satisfies (GA), we obtain 
		$$nT(a)<T(a)^2<e(T(a))=T(a+1)\leq T(b).$$
	\end{proof}

	\begin{remark}\thlabel{rmk:notgrow}
		In general, the conclusion of \Autoref{lem:valpres} does not hold under weaker conditions like $b-a\geq \tfrac 12$. For instance, for a given $a\in \bA^{>0}$ and $c_a\in \varphi(a)$ in \Autoref{cons:trans}, we can set 
		$$f_a\colon [0,1]\to [c_a,e(c_a)]_K, t\mapsto \begin{cases}
		c_a+2t&\!\text{if }0\leq t\leq \tfrac 12,\\
		2(1-t)(c_a+1)+2(t-\tfrac 12)e(c_a)&\!\text{if }\tfrac 12 < t\leq 1.
		\end{cases}$$ Then $f$ maps $[0,\tfrac 12]_K$ bijectively onto $[c_a,c_a+1]_K$ and $[\tfrac 12,1]_K$ onto $[c_a+1,e(c_a)]_K$. The transexponential $T$ thus obtained has the property
		$$T(a)=e(f_a(0))=c_a\text{ and }T(a+\tfrac 12)=e(f_a(\tfrac 12))=c_a+1.$$
		Note that $T$ is not differentiable in the point $x=c_a+1\in \PP_K$.
		\qed 
	\end{remark}

	Recall that we denote by $L=T^{-1}$ the cislogarithm corresponding to $T$. Recall further that $v(\bA^{>0})=G^{<0}$. 

	\begin{definition}
		For any $a\in \bA^{>0}$, we set $$X_T(v(a)):= v(L(a))\in G^{<0}.$$
	\end{definition}
	
	\begin{proposition}\thlabel{prop:xt}
		The map 
		$$X_T\colon G^{<0}\to G^{<0}, g\mapsto X_T(g)$$
		is a contraction on $G^{<0}$. 
	\end{proposition}
	
	\begin{proof}
		Let $a,b\in \bA^{>0}$ with $v(a)=v(b)$. Assume, for a contradiction, that $v(L(a))>v(L(b))$. Since $L(a),L(b)\in \PP_K$, we obtain $L(a)<L(b)$ and $L(b)-L(a)>1$. 
		Hence by \Autoref{lem:valpres}, $v(b)>v(a)$, a contradiction. This shows that $X_T$ is well-defined. In order to show that $X_T$ is surjective, set $c=T(a)\in \PP_K$. Then
		$X_T(v(c))=v(L(c))=v(a)$, as required.
		
		Now let $a,b\in \bA^{>0}$ and set $g=v(a)$ and $h=v(b)$. Suppose that $$X_T(g)=v(L(a))> v(L(b))=X_T(h).$$ Then $L(a)<L(b)$ and thus $a<b$. It now suffices to note that $g=v(a)\geq v(b)=h$. This shows that $X_T$ is monotonically increasing.
		
		We now verify that for any $g,h\in G^{<0}$ with $v_G(g)=v_G(h)$, we have $X_T(g)=X_T(h)$. We may assume that $h<g$. Let $a,b\in \bA^{>0}$ with $v(a)=g$ and $v(b)=h$. Moreover, let $n\in \N^{>1}$ with $ng<h$. 
		First note that $T(nL(a))>T(L(a)+1)=e(T(L(a)))=e(a)>a^n$.  Thus, $nL(a)>L(a^n)$. Since $nv(a)<v(a)<0$, we obtain
		\begin{align*}
		X_T(nv(a))&=X_T(v(a^n))=v(L(a^n))\geq v(nL(a))\\
		&=v(L(a))=X_T(v(a))\geq X_T(nv(a)).
		\end{align*}
		Thus, $X_T(v(a))=X_T(nv(a))$. This yields $X_T(g)=X_T(ng)\leq X_T(h)\leq X_T(g)$ and thus $X_T(g)=X_T(h)$, as required.
	\end{proof}

		Due to \Autoref{prop:xt}, we can now consider the equivalence relation $\sim$ on $G^{<0}$ whose equivalence classes are given by the preimages of $X_T$, that is, for any $g,h\in G^{<0}$ we have \begin{align*}g\sim h :\Leftrightarrow X_T(g)=X_T(h).\end{align*}
		Then $G^{<0}/\!\sim$ is an ordered set with the ordering given by $[g]_\sim<[h]_\sim$ if $X_T(g)<X_T(h)$ for $g,h\in G^{<0}$. Now $G^{<0}$ is isomorphic as an ordered set to $G^{<0}/\!\sim$ via the map
		\begin{align}C_T\colon G^{<0}\to G^{<0}/\!\sim, g\mapsto X_T^{-1}(g).\label{eq:equiv}\end{align}
		This poses a condition on $G$ in order that $(K,e)$ carries a transexponential compatible with $e$. In the following, we will describe $G^{<0}/\!\sim$ more explicitly and also study other necessary conditions on $G$ in order that $(K,e)$ admits a transexponential.
		
		First recall that the map $\bA\to G, a\mapsto -v(e(a))$ is an order-preserving isomorphism of groups.
		This implies that $$\alpha_e\colon \bA^{>0}\to G^{>0}, a\mapsto -v(e(a))$$ defines an order-preserving bijection. 
		As established in \Autoref{prop:condextransexp}, the map
		$$\varphi_T\colon \bA^{>0}\to \PP_K/\!\sim_e, a\mapsto [T(a)]_e$$ defines an order-preserving bijection. 
		Hence, the composition
		$$\varphi_T\circ\alpha_e^{-1}\colon G^{>0}\to \PP_K/\!\sim_e$$
		is an order-preserving bijection between the positive cone of the value group of $K$ and the principal exponential rank of $(K,e)$.
		Now recall that we denote by 
		$$\widetilde{h}_e\colon \Gamma \to G^{<0}, v_G(g)\mapsto v(h_e(g))$$ the group exponential induced by $e$, by
		$$\chi_e\colon G^{<0}\to G^{<0}, v(a)\mapsto v(\ell(a))$$ the corresponding natural contraction %via the group exponential $\widetilde{h}_e\colon \Gamma\to G^{<0}$ 
		and by $$\zeta_e\colon \Gamma\to \Gamma, v_G(g)\mapsto v_G(\chi_e(g))$$ the corresponding shift map. 
		By \cite[Corollary~3.22]{kuhlmann}, the map $$\beta_e\colon\PP_K/\!\sim_e\to G^{<0}/\!\sim_{\chi_e}, [a]_e\mapsto [v(a)]_{\chi_e}$$ is an order-reversing bijection, and the map $$\gamma_e\colon G^{<0}/\!\sim_{\chi_e} \to \Gamma/\!\sim_{\zeta_e}, [g]_{\chi_e}\mapsto [v_G(g)]_{\zeta_e}$$ is an order-preserving bijection.
		We thus obtain the order-preserving maps
		$$\delta_T\colon G^{<0}\to G^{<0}/\!\sim_{\chi_e}, g\mapsto (\beta_e\circ\varphi_T\circ\alpha_e^{-1})(-g)$$
		as well as
		$$\varepsilon_T\colon \Gamma\to \Gamma/\!\sim_{\zeta_e}, \iota \mapsto (\gamma_e\circ\delta_T\circ \widetilde{h}_e)(\iota).$$
		We obtain the following commutative diagram of order-preserving bijections:
		$$
		\begin{tikzcd}
		G^{<0}/\!\sim &G^{<0}  \arrow{r}{\delta_T}\arrow{l}{C_T}
		& G^{<0}/\!\sim_{\chi_e} \arrow{d}{\gamma_e} \\
		&\Gamma \arrow{r}{\varepsilon_T} \arrow{u}{\widetilde{h}_e}
		& \Gamma/\!\sim_{\zeta_e}
		\end{tikzcd}
		$$
		
		\begin{remark}
			In the diagram above, all occuring ordered sets are isomorphic to the principal exponential rank $\PP_K/\!\sim_e$. By abuse of terminology, we also call those sets -- most prominently $G^{<0}/\!\sim_{\chi_e}$ -- the principal exponential rank of $(K,e)$. \qed
		\end{remark}
		
		\begin{lemma}\thlabel{lem:equivcompute}
			Let $a,b\in \bA^{>0}$. Then $v(b)\in \delta_T(v(a))$ if and only if $v(a)=v(e(L(b)))$.
		\end{lemma}
	
		\begin{proof}			
			We first compute $\delta_T^{-1}([v(b)]_{\chi_e})$:
				\begin{align*}
				\delta_T^{-1}([v(b)]_{\chi_e})&= (-\alpha_e\circ\varphi_T^{-1})([b]_e)\\
				&=-\alpha_e(L(b))=v(e(L(b))).
				\end{align*}
				Since $v(b)\in \delta_T(v(a))$ if and only if $[v(b)]_{\chi_e}=\delta_T(v(a))$, by applying $\delta_T^{-1}$ we obtain that $v(b)\in \delta_T(v(a))$ is equivalent to
				$$v(a)=v(e(L(b))).$$
		\end{proof}
	
		Using \Autoref{lem:equivcompute}, we now show that $\sim$ is a coarser equivalence relation than $\sim_{\chi_e}$.
		
		\begin{proposition}
			Let $g,h\in G^{<0}$. Then $g\sim_{\chi_e} h$ implies $g\sim h$. Thus, $[g]_{\chi_e}\subseteq [g]_\sim$.
		\end{proposition}
	
		\begin{proof}
			Let $a,b,c\in \bA^{>0}$ with $v(b)=g$, $v(c)=h$ and $v(b),v(c)\in \delta_T(v(a))$. 
			\Autoref{lem:equivcompute} implies $v(a)=v(e(L(b)))=v(e(L(c)))$.
			Using \Autoref{lem:evgrowth}, we obtain $v(L(b))=v(L(c))$. Thus, $X_T(v(b))=X_T(v(c))$, as required.
		\end{proof}
	
		We now extend \Autoref{prop:condextransexp}.
		
		\begin{theorem}\thlabel{prop:condextransexp2}
			Let $(K,e)$ be a non-archimedean ordered (GAT$_1$)-exponential field. Then the following are equivalent:
			\begin{enumerate}[label = (\roman*)]
				\item 	$(K,e)$ admits a transexponential.
				\item $G^{<0}$ is isomorphic as an ordered set to $G^{<0}/\!\sim_{\chi_e}$.
				\item $\Gamma$ is isomorphic as an ordered set to $\Gamma/\!\sim_{\zeta_e}$.
			\end{enumerate}			
		\end{theorem}
	
		\begin{proof}
			By the discussion above, if $(K,e)$ admits a transexponential, then $\delta_T$ and $\varepsilon_T$ are the required order-preserving bijections. On the other hand, if 
			$\delta\colon G^{<0}\to G^{<0}/\!\sim{\chi_e}$ is an order-preserving bijection, then so is 
			$$\bA^{>0}\to \PP_K/\!\sim_e, a\mapsto \beta_e^{-1}(\delta(-\alpha_e(a))).$$
			This implies that $(K,e)$ admits a transexponential by \Autoref{prop:condextransexp}.
			Finally, if $\varepsilon\colon \Gamma\to \Gamma/\!\sim_{\zeta_e}$ is an order-preserving bijection, then so is 
			$$G^{<0}\to G^{<0}/\!\sim_{\chi_e}, g\mapsto \gamma_e^{-1}(\varepsilon(\widetilde{h}_e^{-1}(g))),$$
			and we are done as above.
		\end{proof}
	
	As a final result in this section, we encode in \Autoref{prop:growthcode} below the growth property of $T$ via $X_T$. 
	
	\begin{lemma}\thlabel{lem:growthtl}
		The following are equivalent:
		\begin{enumerate}[label = (\roman*)]
			\item\label{lem:growthtl:1} $T$ has the growth property.
			
			\item\label{lem:growthtl:2} For any $n\in \N$ and any $a\in \bA^{>0}$, the left-transexponential $\TL$ satisfies $\TL(a)>e^n(a)$.
			
			\item\label{lem:growthtl:3} For any $n\in \N$ and any $a\in \bA^{>0}$ we have $\ell^n(a)>L(a)$.
		\end{enumerate}
	\end{lemma}

	\begin{proof}
		We first show that \autoref{lem:growthtl:1} and  \autoref{lem:growthtl:2} are equivalent.
		\Autoref{prop:rvgrowth} implies that $T$ has the growth property if and only if $\TL$ satisfies $\TL(a)>e^n(a)$ for any $n\in \N$ and any $a\in \PP_K$.
		Suppose that  $\TL(a)>e^n(a)$ for any $n\in \N$ and any $a\in \bA^{>0}$. Let $k\in \N$, let $a\in \bA^{>0}$ and let $b\in R_v$. Then for $n\in \N$ with $n\geq k-\ip{b}+1$ we obtain
		\begin{align*}
			\TL(a+b)&=e^{\ip{b}}(\TL(a+b-\ip{b}))\geq e^{\ip{b}}(\TL(a))>e^{\ip{b}+n}(a)\\
			&\geq e^{k+1}(a)=e^k(e(a))>e^k(2a)>e^k(a+b),
		\end{align*}
		as required.
		
		In order to show that \autoref{lem:growthtl:1} and  \autoref{lem:growthtl:3} are equivalent, we can argue in a similar manner. By setting $c=e^n(a)$ in the above, we obtain that $T$ has the growth property if and only if $\ell^n(c)>L(c)$ for any $n\in \N$ and any $c\in \PP_K$. This shows that \autoref{lem:growthtl:1} implies \autoref{lem:growthtl:3}. The converse can now be shown by a similar computation as above.
	\end{proof}
	
	\begin{proposition}\thlabel{prop:growthcode}
		$T$ has the growth property if and only if for any $n\in \N$ and any $g\in G^{<0}$ we have
		$$X_T(g)>\chi_e^n(g).$$
	\end{proposition}

	\begin{proof}	
		Let $n\in \N$, let $g\in G^{<0}$ and let $a\in \bA^{>0}$ such that $v(a)=g$. If $v(L(a))=X_T(g)>\chi_e^n(g)=v(\ell^n(a))$, then already $L(a)<\ell^n(a)$. Conversely, if $L(a)<\ell^n(a)$, then $X_T(g)\geq\chi_e^n(g)>\chi_e^{n-1}(g)$, as $\chi_e$ is centripetal.
		The proof now follows from \Autoref{lem:growthtl}. 
	\end{proof}

	\Autoref{prop:growthcode} can be rephrased within the principal exponential rank $G^{<0}/\!\sim_{\chi_e}$ of $(K,e)$ as follows.

	\begin{corollary}
		$T$ has the growth property if and only if for any $g\in G^{<0}$ we have $[X_T(g)]_{\chi_e}>[g]_{\chi_e}$.
	\end{corollary}

	\begin{proof}
		Let $g\in G^{<0}$.
		If $[X_T(g)]_{\chi_e}>[g]_{\chi_e}$, the $X_T(g)>\chi_e^n(g)$ for any $n\in \N$. Conversely, if for any $n\in \N$ we have
		$X_T(g)>\chi_e^n(g)$, then $X_T(g)>\chi_e^k(g)$ for any $k\in \Z$, as $\chi_e$ is centripetal. Hence, $X_T(g)\not\sim_{\chi_e}g$, yielding $[X_T(g)]_{\chi_e}>[g]_{\chi_e}$.
	\end{proof}

	%Intoduce here transexponential contraction groups?
	
	\section{Countable models}\label{sec:countable}
	
	We now turn to the study of countable non-archimedean ordered transexponential fields. As a direct consequence of \Autoref{cor:dlo}, we obtain the following.
	
	\begin{proposition}\thlabel{prop:ctbldlo}
		Let $(K,e)$ be a countable non-archimedean ordered \linebreak (GAT$_1$)-exponential field. Then $(K,e)$ admits a transexponential if and only if its principal exponential rank $\PP_K/\!\sim_e$ is isomorphic as an ordered set to $\Q$.
	\end{proposition}

	\Autoref{prop:ctbldlo} raises the question what order types the principal exponential ranks of a countable non-archimedean ordered (GAT$_1$)-exponential fields admits. In the following construction, we show that a countable non-archimedean field admitting one (GAT$_1$)-exponential also admits other (GAT$_1$)-exponentials whose principal exponential ranks have an arbitrary prescribed countable order type.

	\begin{construction}\thlabel{con:exprank}
	Let $K$ be a countable non-archimedean ordered field such that $(K^{>0},\cdot,1,<)$ is divisible and $\ol{K}$ admits a (GAT$_1$)-exponential. Suppose further that $G=\coprod_{\Q}\ol{K}$. By \Autoref{fact:ccc}, $K$ admits a (GAT$_1$)-exponential $\widetilde{e}$. 
	
	Let $\Delta$ be a countable linearly ordered set. 
	We will find a centripetal contraction map $\chi\colon G^{<0}\to G^{<0}$ such that $G^{<0}/\!\sim_\chi$ is isomorphic to $\Delta$ as an ordered set. Before we do so, we first describe how $\chi$ is then used to find a suitable exponential $e$ on $K$ such that the principal exponential rank of $(K,e)$ is isomorphic to $\Delta$.
	
	So suppose that a centripetal contraction map $\chi$ with the property described above is already constructed.
	We set $$h_\chi\colon \Gamma\to G^{<0}, v_G(g)\mapsto \chi(g)$$ for any $g\in G^{<0}$. Then $h_\chi$ is a strong group exponential (see \cite[Section~2.7]{kuhlmann} and \cite[Section~6]{krapp} for details). By \cite[Lemma~4.10]{krapp}, there is a left-exponential $e_{\mathrm{L}}\colon \bA\to \widetilde{e}(\bA)$, that is, an order-preserving isomorphism from the additive ordered group $\bA$ to the multiplicative ordered group $\widetilde{e}(\bA)$, such that its induced group exponential is $h_\chi$ and thus its induced natural contraction is $\chi$. 
	We set $$e(a+b):=e_{\mathrm{L}}(a)\widetilde{e}(b)$$ 	for any $a\in \bA$ and any $b\in R_v$.
	Since $h_\chi$ is a strong group exponential, we obtain that $e$ satisfies (GA) (cf.\ \cite[Section~2.2]{kuhlmann}). Moreover, $e$ satisfies (T$_1$), since its restriction to the valuation ring $e|_{R_v}=\widetilde{e}|_{R_v}$ does so. In conclusion, the principal exponential rank of $(K,e)$ is
	$$G^{<0}/\!\sim_{\chi_e}=G^{<0}/\!\sim_{\chi}\cong \Delta,$$ as required.
	(See also \cite[Sections~1.4~\&~2.5]{kuhlmann} for further details on the construction of exponentials using left-, middle- and right-exponentials.)
	
	We now contruct $\chi$. Let $(\Delta\times \Q,<)$ be ordered lexicographically. Then this is a countable dense linear order without endpoints and we can thus fix an order-preserving bijection $$\eta\colon \Delta\times \Q\to \Q.$$
	For each $\delta\in \Delta$, let $C_\delta\subseteq \Q$ be the convex subset of $\Q$ given by $C_\delta = \eta(\{\delta\}\times \Q)$, and fix some $c_\delta\in C_\delta$. Note that the family $(C_\delta)_{\delta\in \Delta}$ partitions $\Q$ into a family of convex subsets. 
	Now let $\theta\colon \Q\to \Q$ be an order-preserving bijection obtained as follows:
	For any $\delta\in \Delta$, we let $\theta_\delta:=\theta|_{C_\delta}\colon C_\delta\to C_\delta$ be an order-preserving bijection with the property $\theta(q)>q$ for any $q\in C_\delta$ and such that $C_\delta/\!\sim_{\theta_{\delta}}=\{[c_\delta]_{\theta_{\delta}}\}$, i.e.\ $\{\theta_\delta^k(c_\delta)\mid k\in \Z\}$ is coinitial and cofinal in $C_\delta$. Note that this is possible, as $C_\delta$ itself is a dense linear order without endpoints, thus isomorphic to $\Q$, and such a map exists on $\Q$ (e.g. $q\mapsto q+1$).
	We obtain that $\Q/\!\sim_\theta\ =\{[c_\delta]_{\theta}\mid \delta\in \Delta\}$ is isomorphic to $\Delta$ as an ordered set.
	
	Now we construct $\chi\colon G^{<0}\to G^{<0}$ with $G^{<0}/\!\sim_\chi$ isomorphic to $\Q/\!\sim_\theta$. For any $\delta\in \Delta$ and any $k\in \Z$, let $\lambda\colon \Q\to G^{<0}$ be an order-preserving bijection mapping $$[\theta^{k}(c_\delta),\theta^{k+1}(c_\delta)]_\Q\text{ to } [-\one_{\theta^{k+2}(c_\delta)},-\one_{\theta^{k+3}(c_\delta)}]_G.$$ This is possible, as both $\Q$ and $G$ are dense linear orders and the family of convex subsets $([-\one_{\theta^{k+2}(c_\delta)},-\one_{\theta^{k+3}(c_\delta)})_G)_{\delta\in \Delta,\ k\in \Z}$ partitions $G^{<0}$. 
	Finally, we define $\chi\colon G^{<0}\to G^{<0}$ as follows:  
	$$\chi\colon s(q_0)\one_{q_0}+\sum_{q>q_0}s(q)\one_q\mapsto \lambda(q_0)$$
	for any $s(q_0)\in \ol{K}^{<0}$.
	
	We verify that $\chi$ is a centripetal contraction.
	First note that $\chi$ is surjective, as for any $g\in G^{<0}$ we have
	$\chi(-\one_{\lambda^{-1}(g)})=g$. 
	Let $g,h\in G^{<0}$. If $v_G(g)=v_G(h)=q_0\in \Q$, then $\chi(g)=\lambda(q_0)=\chi(h)$. Now suppose that $g<h$ and $q_0=v_G(g)\neq v_G(h)=r_0$. Then $q_0<r_0$. Hence,
	 $\chi(g)=\lambda(q_0)<\lambda(r_0)=\chi(h)$.
	We have thus verified that $\chi$ is a contraction.
	Now let $\delta\in \Delta$ and $k\in \Z$ such that $v_G(g)=q_0\in [\theta^{k}(c_\delta),\theta^{k+1}(c_\delta)]_\Q$. Then $\chi(g)=\lambda(q_0)\in [-\one_{\theta^{k+2}(c_\delta)},-\one_{\theta^{k+3}(c_\delta)}]_G$. Hence, $$v_G(g)\leq \theta^{k+1}(c_\delta)<\theta^{k+2}(c_\delta)\leq v_G(\chi(g)).$$
	This implies $g<\chi(g)$, showing that $\chi$ is centripetal. 
	
	In order to compute $G^{<0}/\!\sim_\chi$, we consider the shift map
	$$\zeta_\chi\colon \Q\to \Q, v_G(g)\mapsto v_G(\chi(g))$$
	for any $g\in G^{<0}$. 
	Following the proof of \cite[Lemma~3.21]{kuhlmann}, we obtain that $G^{<0}/\!\sim_\chi$ and $\Q/\!\sim_{\zeta_\chi}$ are isomorphic as ordered sets. 
	It remains to show that for any $q,r\in \Q$ we have $q\sim_\theta r$ if and only if $q \sim_{\zeta_\chi} r$.
	First suppose that $q \sim_{\zeta_\chi} r$ and assume that $q<r<\zeta_\chi^k(q)$ for some $k\in \N$. Let $\delta\in \Q$ with $q\in [c_\delta]_\theta$. Then
	$$\zeta_\chi(q)=v_G(\chi(-\one_q))=v_G(\lambda(q))\in [c_\delta]_\theta.$$
	Applying this argument iteratively, we obtain $\zeta_\chi^k(q)\in [c_\delta]_\theta$. This implies $r\in [c_\delta]_\theta$, as required. Conversely, suppose that $q\sim_\theta r$ and let $\delta\in \Delta$ with $q,r\in [c_\delta]_\theta$. We may assume that for some $k\in \N$ we have $$\theta^{-k}(c_\delta)<q<r<\theta^k(c_\delta).$$
	We compute $$\zeta_\chi(\theta^{-k}(c_\delta))=v_G(\lambda(\theta^{-k}(c_\delta)))=\theta^{-k+2}(c_\delta)>\theta^{-k+1}(c_\delta)$$
	and obtain iteratively
	$$\zeta_\chi^{2k}(q)>\zeta_\chi^{2k}(\theta^{-k}(c_\delta))>\theta^{k}(c_\delta)>r.$$
	This shows $q\sim_{\zeta_\chi}r$, as required.
	\qed
	\end{construction}

	\begin{remark}\thlabel{rmk:constr}
		\Autoref{con:exprank} can now be used to construct ordered\linebreak (GAT$_1$)-exponential fields $(K,e)$ that admit a transexponential compatible with $e$ and such that do not. For instance, if the construction is started with a finite set $\Delta$, then the resulting $(K,e)$ has finite principal exponential rank and thus does not admit a transexponential by \Autoref{prop:ctbldlo}.
	\end{remark}

	Contrarily, we have the following.

	\begin{lemma}\thlabel{lem:e1per}
		Let $(K,e)$ be a countable non-archimedean ordered (GAT$_1$)-exponential field. Then $K$ admits a (GAT$_1$)-exponential $e_1$ such that the principal exponential rank $\PP_K/\!\sim_{e_1}$ is isomorphic to $\Q$.
	\end{lemma}

	\begin{proof}
		By \Autoref{fact:ccc}, we may assume that $G= \coprod_{\Q} \ol{K}$. Moreover, $(K^{>0},\allowbreak \cdot,1,<)$ is divisible, as this holds for any ordered exponential field, and $\ol{K}$ admits a (GAT$_1$)-exponential $\ol{e}$. Applying \Autoref{con:exprank} to $\Delta=\Q$ gives us the required result.
	\end{proof}

		We now use \Autoref{lem:e1per} to construct an ordered transexponential field $(K,e,T)$ such that $T$ does not have the growth property.
	
	\begin{example}\thlabel{ex:nogrowth}
		Let $(F,\ol{e})$ be a countable archimedean model of $\Th(\R,\exp)$. For instance, $F$ can be the exponential algebraic closure of $\Q$ in $(\R,\exp)$, see \cite{krappprime} for details. Now set $G=\coprod_{\Q}F$ and let $K$ be the real closure of the countable field $F(t^g\mid g\in G)$ in the real closed power series field $F((G))$ (see \cite[Section~1.5]{kuhlmann} for further details on power series fields). Then $\ol{K}=F$, $v(K)=G$, and by \Autoref{fact:ccc}, $K$ admits a (GAT$_1$)-exponential $e$ with residue exponential $\ol{e}$. By \Autoref{lem:e1per}, we can choose $e$ such that the principal exponential rank $\PP_K/\!\sim_e$ is a dense linear order without endpoints. 
		
		We now apply 
		\Autoref{cons:trans} to a suitable order-preserving bijection $\varphi\colon \bA^{>0}\to \PP_K/\!\sim_e$ to obtain a transexponential $T$ without the growth property. For instance, we can pick some $a\in \bA^{>0}$, set $\varphi(a)=[a]_e\in \PP_K/\!\sim_e$ and extend $\varphi$ to an order-preserving bijection from $\bA^{>0}$ to $\PP_K/\!\sim_e$, as both ordered sets are dense linear orders without endpoints.
		
		Finally, using the notation in 	\Autoref{cons:trans}, we have 
		$$T(a)=f_a(0)=c_a\in [a]_e.$$
		Hence, there is some $n\in \N$ with $T(a)<e^n(a)$, but $a\geq (n+1)^2$ as $a$ is a positive infinite element of $K$. This shows that $T$ does not have the growth property.
		\qed
	\end{example}

	In contrast to \Autoref{ex:nogrowth}, we also have that any countable non-ar\-chi\-me\-de\-an ordered transexponential field can be equipped with a possibly different transexponential that has the growth property. This is made precise in the following.

	\begin{proposition}\thlabel{prop:t1growth}
		Let $(K,e,\widetilde{T})$ be a countable non-archimedean ordered transexponential field. Then there exists a transexponential $T$ on $(K,e)$ with the growth property.
	\end{proposition}

	\begin{proof}
		We construct a suitable left-transexponential $\TL\colon \PP_K\to \PP_K$ and then set
		$T$ to be the transexponential on $(K,e)$ given in \Autoref{rmk:alltransexp}. 
		By \Autoref{lem:growthtl}, $\TL$ has to satisfy $\TL(a)>e^n(a)$ for any $n\in \N$ and any $a\in \bA^{>0}$. Constructing $\TL$ such that $[a]_e<[\TL(a)]_e$ for any $a\in \bA^{>0}$ thus suffices. 
		Moreover, using \Autoref{cons:trans}, this reduces to finding a map $\varphi\colon \bA^{>0}\to \PP_K/\!\sim_e$ such that $\varphi(a)>[a]_e$ for any $a\in \bA^{>0}$.
		
		Note that $\PP_K/\!\sim_e$ is a dense linear order without endpoints by \Autoref{prop:ctbldlo}.	
		Let $\eta\colon \Q \to \bA^{>0} $ be an order-preserving bijection. For each $k\in \Z$, let $a_k=\eta(k)$. Now let $b_0\in \PP_K$ with $[b_0]_e>[a_1]_e$. If $b_k$ is already defined for some $k\in \N\cup\{0\}$, let $b_{k+1}\in \PP_K$ with $[b_{k+1}]_e>[b_{k}]_e$ and $[b_{k+1}]_e>[a_{k+2}]_e$. This is possible, as $\PP_K/\!\sim_e$ is a dense linear order without endpoints. Likewise, if  $b_{k}$ is already defined for some $k\in \Z\setminus \N$, let $b_{k-1}\in \PP_K$ with $[b_{k}]_e>[b_{k-1}]_e>[a_k]_e$. We now set
		$$\varphi(a_k)=[b_k]_e$$ for any $k\in \Z$ and 
		extend $\varphi$ to an order-preserving bijection from $\bA^{>0}$ to $\PP_K/\!\sim_e$.
		Now let $a\in \bA^{>0}$ and let $k\in \Z$ with $a_k\leq a<a_{k+1}$. Then
		$$\varphi(a)\geq \varphi(a_k)=[b_k]_e>[a_{k+1}]_e\geq [a]_e,$$
		as required.
	\end{proof}

	We obtain the following characterisation of countable non-ar\-chi\-me\-dean models of $\Totf$.

	\begin{theorem}\thlabel{thm:nonarchexpand}
		Let $K$ be a countable non-archimedean ordered field. Then the following are equivalent:
		\begin{enumerate}[label = (\roman*)]
			\item $\ol{K}$ is exponentially closed in $(\R,\exp)$ and $G=\coprod_{\Q} \ol{K}$.
			
			\item $K$ admits a (GAT$_1$)-exponential.
			
			\item $K$ can be expanded to a model of $T_{\mathrm{otf}}$.
			
			\item  $K$ can be expanded to a model $(K,e,T)$ of $T_{\mathrm{otf}}$ such that $T$ has the growth property.
		\end{enumerate} 
	\end{theorem}

	\begin{proof}
		The equivalences directly follow from \Autoref{lem:gatgat}, \Autoref{fact:ccc}, \Autoref{prop:ctbldlo}, \Autoref{lem:e1per} and \Autoref{prop:t1growth}.
	\end{proof}

	As a final result of this paper, we apply the results of this section to countable non-archimedean models of $\Th(\R,\exp)$. 

	\begin{proposition}
		Let $(K,\widetilde{e})\models \Th(\R,\exp)$ be countable and non-ar\-chi\-me\-de\-an. Then $K$ admits an exponential $e$ and a compatible transexponential $T$ with the growth property such that $(K,e)\models \Th(\R,\exp)$.
	\end{proposition}

	\begin{proof}
		Following \Autoref{con:exprank}, we can find a left-exponential $e_{\mathrm{L}}$ on $K$ such that we obtain an exponential $e$ on $K$ given by $e(a+b)=e_{\mathrm{L}}(a)\widetilde{e}(b)$ for any $a\in \bA$ and $b\in R_v$ and the principal exponential rank of $(K,e)$ is a dense linear order without endpoints. By \Autoref{thm:nonarchexpand} and \Autoref{prop:t1growth}, we can endow $(K,e)$ with a transexponential $T$ that has the growth property. 
		Finally, by Ressayre's axiomatisation \cite{ressayre}, $(K,e)$ is a model of $\Th(\R,\exp)$. (See \cite[Example~5.9]{carlkrapp} for a similar construction and \cite[Section~5.1]{krappthesis} for further details on Ressayre's axiomatisation.)
	\end{proof}

	This section exploited construction methods for ordered transexponential fields starting from countable groups with contractions. Another approach using exponential-logarithmic power series (see, for instance, \cite{kuhlmannkuhlmann}, \cite{kuhlmannshelah} and \cite{kuhlmannmatusinskimantova}) will be subject of a future publication \cite{krappkuhlmannwip}.

		\begin{footnotesize}
		
	\end{footnotesize}


\begin{thebibliography}{99}	
			
			\bibitem{kuhlmannmatusinskimantova}
			\textsc{A.~Berarducci}, \textsc{S.~Kuhlmann}, \textsc{V.~Mantova} and \textsc{M.~Matusinski},
			`Exponential fields and Conway's omega-map',
			\textsl{Proc.\ Amer.\ Math.\ Soc.} 151 (2023) 2655--2669,
			doi:10.1090/proc/14577.
			
			\bibitem{boshernitzan}
			\textsc{M.~Boshernitzan},
			`Hardy fields and existence of transexponential functions',
			\textsl{Aequationes Math.} 30 (1986) 258--280,
			doi:10.1007/BF02189932.
			
			\bibitem{carlkrapp} 
			\textsc{M.~Carl} and \textsc{L.~S.~Krapp}, 
			`Models of true arithmetic are integer parts of models of real exponentiation', 
			\textsl{J. Log. Anal.} 13:3 (2021) 1--21, 
			doi:10.4115/jla.2021.13.3.
				
			\bibitem{dries2} 
			\textsc{L.~van~den~Dries}, \textsc{A.~Macintyre} and \textsc{D.~Marker}, 
			`The elementary theory of restricted analytic fields with exponentiation', 
			\textsl{Ann. of Math. (2)} 140 (1994) 183--205, 
			doi:10.2307/2118545.
			
			\bibitem{krapp}
			\textsc{L.~S.~Krapp},
			`Value Groups and Residue Fields of Models of Real Exponentiation',
			\textsl{J.\ Log.\ Anal.} 11:1 (2019) 1--23,
			doi:10.4115/jla.2019.11.1.
			
			\bibitem{krappthesis}
			\textsc{L.~S.~Krapp},
			\textsl{Algebraic and Model Theoretic Properties of O-minimal Exponential Fields}, Doctoral Thesis, Universität Konstanz, 2019.
			
			\bibitem{krappprime}
			\textsc{L.~S.~Krapp},
			`Embedding the prime model of real exponentiation into o-minimal exponential fields',
			preprint, 2023, 
			arXiv:2302.01609v1.
			
			\bibitem{krappkuhlmannwip}
			\textsc{L.~S.~Krapp} and \textsc{S.~Kuhlmann},
			`Transexponential-cislogarithmic power series fields', in preparation.
			
			\bibitem{kuhlmannkuhlmann}
			\textsc{F.-V.~Kuhlmann} and \textsc{S.~Kuhlmann},
			`Explicit construction of exponential-logarithmic power series'
			\textsl{Proceedings of the Séminaire de Structures Algébriques Ordonnées 1995--1996} (1997) 7 pages.
						
			\bibitem{kuhlmann}
			\textsc{S.~Kuhlmann},
			\textsl{Ordered Exponential Fields},
			Fields Inst.\ Monogr.\ 12 (Amer.\ Math.\ Soc., Providence, RI, 2000),
			doi:10.1090/fim/012.
			
			\bibitem{kuhlmannshelah}
			\textsc{S.~Kuhlmann} and \textsc{S.~Shelah},
			`$\kappa$-bounded exponential-logarithmic power series fields',
			\textsl{Ann.\ Pure Appl.\ Logic} 136 (2005) 284--296,
			doi:10.1016/j.apal.2005.04.001.
			
			\bibitem{miller2} 
			\textsc{C.~Miller}, 
			`Expansions of the real field with power functions', 
			\textsl{Ann.\ Pure Appl.\ Logic} 68 (1994) 79--94, 
			doi:10.1016/0168-0072(94)90048-5.
			
			\bibitem{miller3} 
			\textsc{C.~Miller}, 
			`Exponentiation is hard to avoid', 
			\textsl{Proc.\ Am.\ Math.\ Soc.} 122 (1994) 257--259, 
			doi:10.2307/2160869.
			
			\bibitem{miller}
			\textsc{C.~Miller},
			`A Growth Dichotomy for O-minimal Expansions of Ordered Fields',
			\textsl{Logic: from Foundations to Applications} (eds W.~Hodges, M.~Hyland, C.~Steinhorn and J.~Truss; Oxford Sci.\ Publ., Oxford Univ.\ Press, New York, 1996) 385--399.
			
			\bibitem{padgett}
			\textsc{A.~L.~Padgett},
			\textsl{Sublogarithmic-Transexponential Series}, Doctoral Thesis, University of California, Berkeley, 2022.
			
			\bibitem{ressayre} 
			\textsc{J.-P.~Ressayre}, 
			`Integer Parts of Real Closed Exponential Fields' (extended abstract), 
			\textsl{Arithmetic, Proof Theory and Computational Complexity},  Oxford Logic Guides 23 (eds P.~Clote and J.~Krajícek; Oxford Sci.
			Publ., Oxford Univ. Press, New York, 1993) 278--288.
								
			\bibitem{tarski} 
			\textsc{A.~Tarski}, 
			\textsl{A decision method for elementary algebra and geometry}
			(1948, revised 1951, 2nd edn, RAND Corporation, Santa Monica, Calif., 1957).
			
			\bibitem{wilkie}  
			\textsc{A.~J.~Wilkie}, 
			`Model completeness results for expansions of the ordered field of real numbers by restricted Pfaffian functions and the exponential function', 
			\textsl{J. Amer. Math. Soc.} 9 (1996) 1051--1094, 
			doi:10.1090/S0894-0347-96-00216-0.
			
		\end{thebibliography}
\end{document}